\documentclass[a4paper,fleqn,annotated]{cas-sc}
\usepackage{calc}
\usepackage{etoolbox}
\usepackage{graphicx}
\usepackage{xcolor}

\makeatletter
\newlength{\hb@colwidth}
\def\makeheaderbox{
	\noindent\textcolor{black}{\rule{\textwidth}{1pt}}
	\par\vskip2pt\noindent
	\begin{minipage}{\textwidth}
		\centering
		\setlength{\tabcolsep}{0pt}
		\setlength{\extrarowheight}{0pt}
		\setlength{\hb@colwidth}{0.16666667\textwidth}
		\renewcommand{\arraystretch}{2}
		\begin{tabular}{
				p{\hb@colwidth}
				>{\columncolor[gray]{0.95}}p{4\hb@colwidth}
				p{\hb@colwidth}
			}
			\parbox[c][\hb@colwidth][c]{\hb@colwidth}{
				\centering
				\includegraphics[height=0.9\hb@colwidth, keepaspectratio]{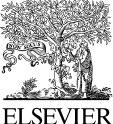}
			} & 
			\parbox[c][\hb@colwidth][c]{4\hb@colwidth}{
				\centering
				{\fontsize{8}{10}\selectfont Contents lists available at \href{https://www.elsevier.com/locate/fi}{ScienceDirect}}\par
				\vspace{15pt}
				{\fontsize{16}{18}\selectfont Journal of the Franklin Institute}\par
				\vspace{12pt}
				{\fontsize{8}{10}\selectfont journal homepage: \href{https://www.elsevier.com/locate/fi}{www.elsevier.com/locate/fi}}
			} & 
			\parbox[c][\hb@colwidth][c]{\hb@colwidth}{
				\centering
				\includegraphics[height=\hb@colwidth, keepaspectratio]{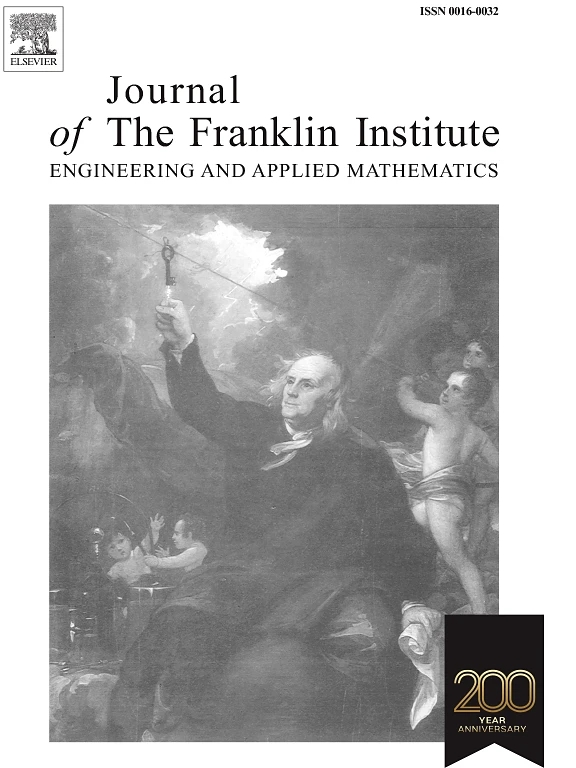}
			} \\
		\end{tabular}
	\end{minipage}
	\par\vskip2pt
	\noindent\textcolor{black}{\rule{\textwidth}{3pt}}
	\par\vskip18pt
}

\let\originalmaketitle\maketitle
\def\maketitle{
	\makeheaderbox
	\vspace{2pt}
	\begin{center}{\Large\bfseries \emph{Bicentennial Special Issue Invited Paper}}
    \end{center}\par
	\vspace{10pt}
	\originalmaketitle
}
\makeatother

\usepackage[numbers]{natbib}
\def\tsc#1{\csdef{#1}{\textsc{\lowercase{#1}}\xspace}}
\tsc{WGM}
\tsc{QE}
\tsc{EP}
\tsc{PMS}
\tsc{BEC}
\tsc{DE}

\usepackage{hyperref}
	\hypersetup{hidelinks=true}
\usepackage{mathtools}
\usepackage{tikz}
	\usetikzlibrary{arrows.meta}

\newtheorem{thm}{Theorem}
\newtheorem{lem}{Lemma}
\newtheorem{prop}{Proposition}

\newtheorem{asmp}{Assumption}
\newtheorem{defn}{Definition}

\newtheorem{rem}{Remark}
\newproof{pf}{Proof}

\DeclareMathOperator{\vecop}{vec}

\DeclareMathOperator{\vech}{vech}

\newcommand{\revcolor}{red!60!black}
\newcommand{\rev}[1]{\textcolor{\revcolor}{#1}}

\begin{document}
\let\WriteBookmarks\relax
\def\floatpagepagefraction{1}
\def\textpagefraction{.001}
\shorttitle{Journal of the Franklin Institute}
\shortauthors{McNamee et~al.}

\title[mode=title]{
Extremum Seeking (ES) is Practically Stable Whenever Model-Based ES is Stable
}
\tnotemark[1]

\author[1,2]{Patrick McNamee}[orcid=0000-0001-9303-6931]

\author[1]{Zahra Nili Ahmadabadi}[orcid=0000-0003-3226-4425]

\author[2]{Miroslav Krsti{\'c}}[orcid=0000-0002-5523-941X]
\cormark[1]

\affiliation[1]{
	organization={Department of Mechanical Engineering, San Diego State University},
	city={San Diego},
	state={California},
	country={United States}
}

\affiliation[2]{
	organization={Department of Mechanical and Aerospace Engineering, University of California San Diego},
	city={La Jolla},
	state={California},
	country={United States},
}

\cortext[cor1]{Corresponding author}
\tnotetext[1]{Research was sponsored by the Army Research Office under Grant
Number W911NF-24-1-0386 and the Department of the Navy, Office of Naval Research
under ONR award number N000142412269. The views, findings, conclusions, or recommendations contained in this document are those of the authors and should not be interpreted as representing the official policies or views, either expressed or implied, of the Army Research Office, the Office of Naval Research, or the U.S. Government. The U.S. Government is authorized to reproduce and distribute reprints for Government purposes notwithstanding any copyright notation herein.}

\begin{abstract}
\rev{Extremum seeking controllers (ESCs) are} optimization algorithms in continuous time, with model-based ESCs using true derivative information of the cost function and model-free ESCs utilizing perturbation-based estimates instead. Stability analysis of model-free ESCs often employs the associated average system, whose stability is dependent on the selection of the dither signal. We first demonstrate the challenge of this analysis approach by showing selections of relative dither amplitudes and rates at different ESC inputs which result in the average system always having an unstable equilibrium. We then show that, under averaging hypotheses that hold for the standard ESC constructions, if the model-based ESC is a continuous (but not necessarily locally Lipschitz) vector field that is globally asymptotically stable (GAS), then the average system is semiglobally practically asymptotically stable (sGPAS), and the model-free ESC is semiglobally practically uniformly asymptotically stable (sGPUAS). Thus, the system analyst does not need to analyze the stability of each average system that corresponds to the choice of dither signals, as the stability of the model-based ESC is sufficient to show that all the average systems are sGPAS for appropriately selected dither signals.
\end{abstract}


\maketitle

\section{Introduction}
	\label{sec:intro}
    
Extremum Seeking Control (ESC) is a family of continuous-time optimization algorithms which have been used for a variety of applications such as: automatic braking systems \cite{ref:zhang-2007}; maximum power point tracking for photovoltaic arrays \cite{ref:leyva-2006}; and stabilization of combustion instability \cite{ref:banaszuk-2000}. These algorithms seek to find extrema, preferably global extremum, of some sensor field or cost function $J$ that depends on a parameter vector $\theta$. We refer to ESCs which can directly obtain knowledge of $J$'s derivatives as model-based ESCs whereas ESCs which do not have such knowledge must rely on perturbation-based estimates and are thus referred to as model-free ESCs. 
	
The first practical stability proof for a model-free ESC was not a global stability proof but rather a local practical stability proof of a gradient-based ESC (GESC) \cite{ref:krstic-2000}. This proof relied on a Taylor series expansion of the map around an extremum to guarantee stability using the quadratic term in the map and estimate the direction and magnitude of the bias using the cubic term. The later work \cite{ref:tan-2006} extended this result from local to semiglobal by making  the dither amplitude and seeking speed sufficiently small. However, these results assumed that the cost function to be minimized satisfies $(\theta - \theta^*)^T \nabla J(\theta) > 0$ for all $\theta \neq \theta^*$ where $\theta^*$ is the global minimizer. It is unclear how to extend these results to other algorithms, such as Newton-based ESC (NESC) \cite{ref:ghaffari-2012,ref:mcnamee-2024, ref:mcnamee-2025}, where there may be additional auxiliary states to consider. Furthermore, the assumptions in \cite{ref:tan-2006} impose restrictions on $J$ which may be unnecessarily restrictive.
	
The articles \cite{ref:nesic-2010,ref:nesic-2013} consider frameworks for generating practically stable model-free ESCs from model-based ESCs. These considered frameworks for generating model-free ESCs have the derivatives of $J$ being estimated on a different, faster time scale than the parameter dynamics but slower than the dither time scale. The motivation was that the stability of the overall model-free ESC could be determined by the stability of both the derivative estimators and the model-based ESCs. However, there are several drawbacks to utilizing a third intermediate time scale. The first drawback is theoretical: it complicates the proof of practical stability, requiring singular perturbation techniques to analyze both a boundary layer model and a reduced order model. This additional analysis imposes unnecessary assumptions on the model-based system to ensure the practical stability of the reduced order model. Moreover, it introduces an unnecessary small parameter which is the root cause of the second and third drawbacks. These next two drawbacks relate to the effort required to tune parameters for ESC for deployment in the real-world and the subsequent real-world performance. The tuning of small parameters for ESCs in the real-world is often done in reverse compared to the tuning of the small parameters in practical stability analyses. For the methodology given in \cite{ref:nesic-2010,ref:nesic-2013}, the system designer selects an initial dither amplitude before executing the following iterative tuning loop:
\begin{enumerate}
    \item First, the designer sets the dither rate based on any hardware/software speed limitations;
    \item Then, the designer ensures that the intermediate time scale is sufficiently slower than the fast time scale so that the boundary layer model has sufficiently small oscillations around the quasi-steady state in the regions of interest;
    \item Next, the designer selects gains so that the slow time scale is sufficiently slower than the intermediate time scale for the reduced order model to be a reasonable approximation to the slow states of the average system;
    \item Lastly, the designer checks if the ESC converges to a satisfactorily small neighborhood of the extremum. If the final convergence neighborhood is too large, the dither amplitude is reduced and tuning starts again from step 2.
\end{enumerate}
The extraneous effort in tuning the small parameters is found in steps 2 and 3: more effort/time is required to tune systems with three time scales than those with only two. Even if this tuning effort is overcome, the third drawback limits real-world performance. Ideally, the system designer wants to tune gains for faster parameter convergence by speeding up the slow time scale while, at the same time, conserving the size of the region of attraction and the asymptotic convergence accuracy. However, this goal is countered by two limitations: the hardware/software limits that place upper bounds on the speed of the fast time scale, and the requirement that the intermediate time scale be strictly between the fast and slow time scales. This relative placement of the intermediate time scale limits the speed of the slow time scale, since strict gaps must be maintained between the slow and intermediate time scales, as well as between the intermediate and fast time scales. Removal of the intermediate time scale allows for the slow time scale to be limited instead by the quicker fast time scale, thus improving the realized convergence speed of the parameters. 

Additionally, these previous framework approaches exclude other optimization methods such as ADAM \cite{ref:kingma-2015}, AdaGrad \cite{ref:duchi-2011}, and RMSprop \cite{ref:ma-2022} which are all simultaneously estimating derivatives and updating parameters on the same time scale. We do not claim that these discrete-time methods are themselves instances of the abstraction given in Section~\ref{sec:extremum-seekers}; what the framework accommodates are suitably regularized continuous-time analogues of them, whose auxiliary states evolve on the same time scale as the parameter estimate and whose vector fields satisfy mild assumptions. One could attempt to analyze these additional systems with averaging theory, but this analysis would likely rely on the average system stability to conclude the practical stability of the model-free ESC. The average system, as we show later, may not be globally asymptotically stable (GAS) for certain dither signals when acting on a multivariable $J$ that is not quadratic. The system analyst would need to determine for each $J$ what dither signals, if any, would result in the average system being stable. This task is made even harder as one does not know $J$ before deploying the model-free ESC.

Rather than utilizing the existing methodologies in the framework approaches, where the start of the system design is a model-based ESC, it is possible for the system designer to start with the construction of a model-free ESC and derive the model-based ESC. This approach eliminates the intermediate time scale, reducing the required number of small parameters for practical stability from three in the framework approach to just two. In fact, as we show in Section~\ref{sec:extremum-seekers:mb-from-mf} using results of \cite{ref:kavanaugh-2025}, it is relatively straightforward to construct model-free ESCs from which the model-based ESCs can be quickly derived without calculation of the average system, although this model-based ESC may not be a previously studied optimization algorithm. If this potentially unstudied model-based ESC is GAS, then the model-free ESC is practically stable with respect to two parameters. This is the main result of this work, and its proof is presented in Section~\ref{sec:main-results}. The proof utilizes what the system analyst has likely found: a Lyapunov function that proves the model-based ESC is GAS. This differs from the proofs in \cite{ref:nesic-2010,ref:nesic-2013}, which use class $\mathcal{KL}$ comparison functions to serve as outer envelopes for the trajectories of the model-based ESC system. We prefer the use of a Lyapunov function because comparison functions tend to be implicitly defined based on differential inequalities rather than explicitly constructed; however, for completeness, we provide $\mathcal{KL}$-based proofs in Appendix~\ref{sec:kl-based-proofs}. Unlike previous literature, the proof of our main result relies on general averaging theory rather than the standard periodic averaging theory, which reduces the restrictions on the designer's choice of dither signals.

\emph{Relation to prior work:} The implication ``GAS of a limiting vector field $\Rightarrow$ practical stability of a fast time-varying realization'' has a long history, and we do not claim the mechanism itself as new. Teel, Peuteman, and Aeyels \cite{ref:teel-1998,ref:teel-1999} established that GAS of the average implies semiglobal practical asymptotic stability of the actual system; Ne{\v s}i{\'c} and Teel \cite{ref:nesic-2001} and Teel and Ne{\v s}i{\'c} \cite{ref:teel-2000} developed averaging and closeness-of-solutions results under strong and weak averages with disturbances; Teel, Moreau, and Ne{\v s}i{\'c} \cite{ref:teel-2003} unified the two-time-scale ISS statements; Moreau and Aeyels \cite{ref:moreau-2000} treated practical stability for converging families of systems; and Tan, Ne{\v s}i{\'c}, and Mareels \cite{ref:tan-2005,ref:tan-2006} specialized such statements to extremum seeking. What we contribute is not a new robustness mechanism but an ESC-specific synthesis with four elements: (i) the decisive hypothesis is identified as \emph{compact} (locally uniform) convergence of the averaged vector field to the model-based one, a property which is a consequence of the demodulation construction of a model-free ESC, so that nothing needs to be verified about the model-based system beyond its stability; (ii) the model-based vector field is required only to be continuous, not locally Lipschitz, so that closeness of trajectories between the average and model-based systems is unavailable and a Lyapunov argument is used instead; (iii) the dither is allowed to be almost periodic rather than periodic, which requires averaging over an infinite interval; and (iv) the proof is written so that the analyst reuses the Lyapunov function already constructed for the model-based system, rather than a comparison function that is only implicitly defined.

\emph{Main Contribution:}
The  main contribution of this work is to eliminate the need for the system analyst to consider every possible dither signal when determining the practical stability of the ESC algorithm. We show that, under conditions outlined in  Assumption~\ref{asmp:averaging} in Section~\ref{sec:extremum-seekers:abstraction}, the model-based ESC being a continuous, but not necessarily locally Lipschitz, GAS vector field is sufficient to guarantee that the model-free ESC is semiglobally practically uniformly asymptotically stable (sGPUAS), even if the extremum may be unstable for an average system with specific dither signals. This result follows from the compact convergence of the average system to the model-based system and general averaging techniques. This extends the GESC semiglobal practical stability results in \cite{ref:tan-2006} to other model-free ESC systems, systems which may be perturbation-based equivalents of existing model-based ESCs systems or completely novel model-free ESCs themselves.

\emph{Notation:} We use $\Vert\cdot\Vert$ to denote the standard Euclidean norm and $\mathcal{B}_{\rho}$ to denote the open ball $\lbrace\ x\ \vert\ \Vert x \Vert < \rho \rbrace$.

\section{Extremum Seekers}
    \label{sec:extremum-seekers}
    
ESC algorithms act on cost functions $J:\mathbb{R}^n\to\mathbb{R}$ that map a parameter $\theta\in\mathbb{R}^n$ to some real number output. The ESC moves some estimate of the optimal parameter $\hat{\theta}$ to the optimal parameter $\theta^*$ which, in this work, minimizes $J$. To accomplish this in a model-free way, one must perturb $\theta$ around $\hat{\theta}$ to gain local knowledge of $J$ to estimate the derivatives, with it being generally assumed that $J$ is sufficiently differentiable so that the derivatives in question exist. To this extent, the perturbed parameter is 
\begin{equation}
    \theta(t) = \hat{\theta}(t) + a s(\omega t)\,,
\end{equation} 
where $s(\omega t)$ is a baseline dither signal whose trajectories are shaped by the system designer. The parameters $a > 0$ and $\omega > 0$ are the small and large parameters, respectively, and control the absolute size and speed of the chosen additive dither signal. The dither signal $s$ is a continuous function and defined element-wise for sinusoidal perturbations as 
\begin{equation}
    a s_i(\omega t) = a r_i\sin(\omega'_i \omega t) = a_i \sin(\omega_i t)\,.
\end{equation} 
Here, $\omega'_i$ and $r_i$ control the relative dither rates and amplitudes while $\omega_i$ and $a_i$ are the absolute dither rates and amplitudes. As in \cite{ref:ghaffari-2012}, the relative dither rates must be restricted  such that $\omega'_i\neq \omega'_j$ and $\vert \omega'_i \vert \neq 2\vert \omega'_j\vert$ for $i\neq j$ although some ESCs may require additional restrictions depending on the derivatives being estimated for the ESC (cf.~\cite{ref:ghaffari-2024}). Furthermore, the designer must choose the relative dither amplitudes $r_i$ such that $\sum_{i=1}^n r_i^2 = 1$ and $r_i\neq 0$ for all $i$ so that $s$ is normalized. However, unlike previous literature, we do not restrict the ratios $\omega_i'/\omega_j'$ to be rational and instead allow them to potentially be irrational. This means that the vector $s$ may be composed of periodic elements, but there may not exist a common period between all elements. The lack of a common period means that $s$ is almost-periodic. Consequently, the averaging theory used in our analysis requires averaging over an infinite time interval, rather than a finite time interval used in previous literature \cite{ref:ghaffari-2012,ref:nesic-2010,ref:nesic-2013}, as well as justification that the vector fields are almost periodic with respect to (wrt) the fast time scale $\tau = \omega t$ whenever $s$ is almost periodic. A benefit of considering almost periodic dither signals and vector fields in our analysis is the generality of our result: periodic functions are a subset of almost periodic functions, so our results cover a broader family of extremum seekers.

While $s$ is commonly constructed with sinusoidal functions in the literature, alternative dither signals may be utilized, provided they are continuous periodic or almost periodic functions \cite{ref:tan-2007,ref:haring-2023,ref:kavanaugh-2025}. In particular, \cite{ref:kavanaugh-2025} gives the necessary and sufficient condition for determining if an $s$ can be used to demodulate the necessary derivative information for any model-free ESC. As noted in \cite{ref:tan-2007} and demonstrated in \cite{ref:kavanaugh-2025}, the choice of $s$ and how derivatives are demodulated will affect the performance of the model-free system, such as the region of attraction of $\theta^*$, the speed of convergence, and the ultimate convergence bound of the model-free ESC trajectories. Given the prevalence of sinusoidal perturbations in previous literature and for ease of reference, we state the restrictions on $\omega_i'$ sufficient for a sinusoidal $s$ in the model-free ESCs shown in Sections \ref{sec:extremum-seekers:gesc} and \ref{sec:extremum-seekers:nesc}. If a system designer wants to use periodic sinusoidal dither signals, then they must find a set of $\omega'_i$ which satisfy these rules. However, if the system designer opts for almost periodic dither signals and employs a sequence provided in \cite{ref:kavanaugh-2025}, such as $\omega'_i = 2^{(i-1)/n}$, then this choice is guaranteed to satisfy the rules given later.

With the perturbations of the parameter defined, we now go through some model-free ESC systems and how they use the cost function measurement to update the estimates of the optimizing input parameters.

\subsection{Gradient-Based ESC}
    \label{sec:extremum-seekers:gesc}

The baseline ESC algorithm considered in this work is the gradient-based ESC (GESC), which is intended to work on $J$'s which are continuously differentiable ($J\in\mathcal{C}^1$). The model-free, perturbation-based GESC updates $\hat{\theta}$ by the following differential equation
\begin{equation}
    \label{eq:gesc:model-free:parameter-update}
	\frac{d}{dt}\hat{\theta} = -k \hat{g}(\omega t, \hat{\theta}, a)
\end{equation}
where $\hat{g}$ is the estimate of the gradient $\nabla J$ formed from the perturbed $J$. It is defined element-wise as
\begin{equation}
	\hat{g}_i(\tau, \hat{\theta}, a) =  m_i(\tau, a)J(\hat{\theta} + a s(\tau))
\end{equation}
with $m_i$ being the elements of the time-varying vector function $m(\tau,a)$ multiplied by the perturbed output of $J$ to form the perturbation-based estimate of the gradient. For sinusoidal $s$, $m$ is commonly defined element-wise as
\begin{equation}
	m_i(\tau, a) = \frac{2}{a r_i} \sin(\omega'_i \tau)\,.
\end{equation}
To apply averaging theory whenever $s$ is almost periodic, we must justify that \eqref{eq:gesc:model-free:parameter-update} is almost periodic wrt $\tau$. Since all $m_i$'s are periodic wrt $\tau$ and, for fixed $\hat{\theta}$, the perturbed $J$ is almost periodic wrt $\tau$ whenever $J$ is continuous, every $\hat{g}_i$ is almost periodic wrt $\tau$ \cite[Th~1.9]{ref:fink-1975}. Consequently, $\hat{g}$ is almost periodic wrt $\tau$ \cite[Th~2.9]{ref:fink-1975}. This means that the associated average system is
\begin{equation}
    \label{eq:gesc:average-system:parameter-update}
	\frac{d}{dt} \hat{\bar{\theta}} = -k \hat{\bar{g}}(\hat{\bar{\theta}}, a)
\end{equation}
where average gradient estimate $\hat{\bar{g}}$ is defined as
\begin{equation}
    \label{eq:gesc:average-gradient-estimate}
	\hat{\bar{g}}(\hat{\bar{\theta}}, a) = \lim_{T\to\infty} \frac{1}{T}\int_{0}^{T} \hat{g}(\tau, \hat{\bar{\theta}}, a) d\tau = \nabla J(\hat{\bar{\theta}}) + R_{g}(\hat{\bar{\theta}}, a)
\end{equation}
with the integral being well-defined by \cite[Th~3.1]{ref:fink-1975} and $R_g$ being the difference between the average gradient estimate and the true gradient ($R_g = \hat{\bar{g}} - \nabla J$). Note that while \eqref{eq:gesc:average-gradient-estimate} is integrated over the interval $[0,T]$, the limit in \eqref{eq:gesc:average-gradient-estimate} converges uniformly when integrating over the interval $[t_0, t_0+T]$ for any $t_0\in\mathbb{R}$ \cite[Cor~3.2]{ref:fink-1975}. Additionally, for the choice of an additive sinusoidal dither and previously mentioned restrictions on $\omega_i'$'s, we know that $R_g(\hat{\bar{\theta}}, a) \to 0$ as $a\to 0$ for all $\hat{\bar{\theta}}\in\mathbb{R}^n$ if $J\in\mathcal{C}^1$ \cite[Th~1]{ref:kavanaugh-2025}. Hence, we define the model-based GESC as the system governed by the dynamics in the limit $a\to 0$, which is the classical gradient-descent differential equation
\begin{equation}
    \label{eq:gesc:model-based:parameter-update}
	\frac{d}{dt}\vartheta = -k \nabla J(\vartheta)\,.
\end{equation}

It is important to note the quality of convergence of the average system in \eqref{eq:gesc:average-system:parameter-update} to the model-based system in \eqref{eq:gesc:model-based:parameter-update}. Clearly, the convergence is dictated by the convergence of $R_g(\hat{\bar{\theta}}, a)\to 0$ as $a\to 0$, but can this be strengthened from just point-wise convergence? It was shown in \cite{ref:kavanaugh-2025} that the limit $R_g(\hat{\bar{\theta}}, a)\to 0$ as $a\to 0$ has, in a worst-case analysis, the same order of convergence as $\nabla J(\theta) \to \nabla J(\hat{\bar{\theta}})$ when $\theta \to \hat{\bar{\theta}}$. Since $J\in\mathcal{C}^1$, $\nabla J$ is a continuous function. Therefore, on any compact set $\mathcal{D}\subset\mathbb{R}^n$, we can always find a function $\alpha\in\mathcal{K}$ such that $\Vert \nabla J(\theta) - \nabla J(\hat{\bar{\theta}}) \Vert \leq \alpha (\Vert \theta - \hat{\bar{\theta}} \Vert)$ for all $\theta, \hat{\bar{\theta}}\in\mathcal{D}$, and a constant $c\geq 1$ such that $\Vert R_g(\hat{\bar{\theta}}, a) \Vert \leq c \alpha(a)$ for all $\hat{\bar{\theta}}\in\mathcal{D}$.

Another critical consideration is the regularity required for the vector fields of these systems. For the existence of the average and model-based systems in \eqref{eq:gesc:average-system:parameter-update} and \eqref{eq:gesc:model-based:parameter-update}, it is necessary to assume that $J$ is differentiable. Differentiability of $J$ is sufficient to conclude that $\hat{g}$ is locally Lipschitz in $\hat{\theta}$ for a fixed $a$ and, therefore, the trajectories of the model-free system are unique. Similarly, $J$ being differentiable implies that $\hat{\bar{g}}$ is also differentiable with respect to $\hat{\bar{\theta}}$ as well as the trajectories of the average system being unique. In contrast, we have assumed little about \eqref{eq:gesc:model-based:parameter-update}, only that the vector field is continuous since $\nabla J$ is continuous. Consequently, we are assuming that the solutions exist for the model-based system but not that these solutions are unique. This lack of uniqueness prevents the establishment of a ``closeness of trajectories'' between trajectories of the average and model-based systems. One cannot ensure that a trajectory of the average system stays sufficiently close, given sufficiently small $a$, to any specified trajectory of the model-based system if multiple trajectories can occur for the same initial condition $\vartheta(t_0) = \vartheta_0$. It turns out that this lack of closeness between the average system and the model-based system is not an issue for us; Theorem~\ref{thm:model-based-implies-model-free} gives the result that if the model-based system is GAS, then the average system is sGPAS and the model-free system is sGPUAS.

\subsection{Newton-Based Extremum Seeking Control}
	\label{sec:extremum-seekers:nesc}

A more complicated ESC is the Newton-based ESC (NESC) \cite{ref:ghaffari-2012}, 
\begin{eqnarray}
	\label{eq:nesc:model-free:parameter-update-ode}
	\frac{d}{dt} \hat{\theta} &=& -k\hat{\Gamma} \hat{g}(\omega t, \hat{\theta}, a) \\
	\label{eq:nesc:model-free:inverse-hessian-ode}
	\frac{d}{dt}\hat{\Gamma} &=& \omega_l \left(\hat{\Gamma} - \hat{\Gamma}\ \hat{H}(\omega t, \hat{\theta}, a)\ \hat{\Gamma}\right)
\end{eqnarray}
where $\hat{\Gamma}$ is the estimate of the matrix inverse of the Hessian estimate $\hat{H}$ of the true unknown Hessian $\nabla^2 J$. With sinusoidal dither signals, the Hessian is estimated element-wise by
\begin{align}
    \label{eq:nesc:model-free:hessian-diagonal}
	\hat{H}_{ii}(\tau, \hat{\theta}, a) &{}={} \frac{16}{a^2 r_i^2}\left(\sin^2(\omega'_i \tau) - \frac{1}{2}\right) J(\hat{\theta} + a s(\tau)) \\
    \label{eq:nesc:model-free:hessian-off-diagonal}
	\hat{H}_{ij}(\tau, \hat{\theta}, a) &{}={} \frac{4}{a^2 r_i r_j}\sin(\omega'_i \tau)\sin(\omega'_j \tau) J(\hat{\theta} + a s(\tau))
\end{align}
and, for logic identical to $\hat{g}$, is almost periodic wrt $\tau$. Since second-order derivatives are being estimated, there need to be additional restrictions on the relative dither rates $\omega'_i$. The additional restrictions on the relative dither rates are that $\vert \omega'_i \vert \pm \vert \omega'_j \vert \neq \vert \omega'_k \vert$, $\vert \omega'_i \vert \pm \vert \omega'_j \vert \neq 2\vert \omega'_k \vert$, and $\vert \omega'_i \vert \pm \vert \omega'_j \vert \neq \vert \omega'_k \vert \pm \vert \omega'_\ell \vert$ for all distinct $i$, $j$, $k$, and $\ell$ \cite{ref:ghaffari-2012}. 

Stability of the NESC is often focused on the average system, which is given by
\begin{eqnarray}
	\label{eq:nesc:average-system:parameter-update-ode}
	\frac{d}{dt} \hat{\bar{\theta}} &=& -k\hat{\bar{\Gamma}} \hat{\bar{g}}(\hat{\bar{\theta}}, a) \\
	\label{eq:nesc:average-system:inverse-hessian-ode}
	\frac{d}{dt}\hat{\bar{\Gamma}} &=& \omega_l \left(\hat{\bar{\Gamma}} - \hat{\bar{\Gamma}}\ \hat{\bar{H}}(\hat{\bar{\theta}}, a)\ \hat{\bar{\Gamma}}\right)
\end{eqnarray}
where the average Hessian estimate $\hat{\bar{H}}$ is defined as
\begin{align}
    \label{eq:nesc:average-hessian-estimate}
	\hat{\bar{H}}(\hat{\bar{\theta}}, a) &{}={} \lim_{T\to\infty} \frac{1}{T} \int_{0}^T \hat{H}(\tau, \hat{\bar{\theta}}, a) d\tau = \nabla^2 J(\hat{\bar{\theta}}) + R_{H}(\hat{\bar{\theta}}, a)
\end{align}
where the limit is again well defined \cite[Th~3.1]{ref:fink-1975} and $R_H$ is the difference between the average Hessian estimate and the true Hessian ($R_H = \hat{\bar{H}} - \nabla^2 J$). For similar reasons to $R_g$, $R_H(\cdot, a)\to 0$ compactly as $a\to 0$ but only if $J$ is twice differentiable ($J\in\mathcal{C}^2$) \cite{ref:ghaffari-2012, ref:kavanaugh-2025}. Thus, the model-based NESC, which is defined by taking the limit of $a\to 0$ of \eqref{eq:nesc:average-system:parameter-update-ode} and \eqref{eq:nesc:average-system:inverse-hessian-ode}, is
\begin{eqnarray}
	\label{eq:nesc:model-based:parameter-update-ode}
	\frac{d}{dt} \vartheta &=& -k \Pi \nabla J(\vartheta) \\
	\label{eq:nesc:model-based:inverse-hessian-ode}
	\frac{d}{dt}\Pi &=& \omega_l \left(\Pi - \Pi\ \nabla^2 J(\vartheta)\ \Pi\right)
\end{eqnarray}

Although the NESC is composed of a vector differential equation and a matrix differential equation, it has a reformulation in which the NESC can be described by two vector differential equations. This reformulation allows for more familiar stability results to be applicable. For the reformulation in this work, we will assume that the NESC is only acting on cost functions with strictly positive definite Hessians ($\nabla^2 J(\hat{\bar{\theta}}) \succ 0$ $\forall\hat{\bar{\theta}}\in\mathbb{R}^n$) and so all trajectories of $\hat{\Gamma}(t)$ will be positive definite symmetric matrices as long as the initial estimate $\hat{\Gamma}(t_0)$ is a positive definite symmetric matrix \cite{ref:ghaffari-2012}. Two qualifications are needed here, since they bear on whether the model-free NESC is an instance of the abstraction of Section~\ref{sec:extremum-seekers:abstraction}. First, positive definiteness of $\hat{\bar{H}}$ for forward complete trajectories is an invariance property of the \emph{average} filter \eqref{eq:nesc:average-system:inverse-hessian-ode}  for sufficiently small $a$. The equivalent invariance is not true for the model-free filter \eqref{eq:nesc:model-based:inverse-hessian-ode} as the instantaneous demodulated estimate $\hat{H}(\omega t, \hat{\theta}, a)$ in \eqref{eq:nesc:model-free:hessian-diagonal}--\eqref{eq:nesc:model-free:hessian-off-diagonal} is sign indefinite for every $J$ due to the demodulating signals having zero mean. Second, positive definiteness of $\hat{\Gamma}$ is nevertheless preserved along \eqref{eq:nesc:model-free:inverse-hessian-ode}, because $v^\top \hat{\Gamma} v = 0$ implies $v^\top \frac{d}{dt}\hat{\Gamma} v = 0$, so the boundary of the positive semidefinite cone is invariant; what is \emph{not} guaranteed is forward completeness. In the scalar case with $\hat{H} \leq -c < 0$ on an interval, \eqref{eq:nesc:model-free:inverse-hessian-ode} gives $\frac{d}{dt}\hat{\Gamma}\geq \omega_l c \hat{\Gamma}^2$ and $\hat{\Gamma}$ escapes in finite time; equivalently, $P = \hat{\Gamma}^{-1}$ obeys $\dot{P} = \omega_l(\hat{H} - P)$ and is driven through zero. This is not a pathological corner: for $J(\theta) = 1 + \theta^2$, trajectories with initial conditions satisfying $\hat{\Gamma}(t_0) \geq \sqrt{2} \omega a^2 /\pi\omega_{\ell}$ have finite escape times within $\pi/4\omega$ seconds since $\hat{H}\leq -4\sqrt{2} a^{-2}$ within this time interval. Accordingly, the NESC serves as an example of the \emph{form} of the abstraction in Section~\ref{sec:extremum-seekers:abstraction}, and the logarithmic coordinates are used on a forward invariant compact set of $\hat{\Gamma}$. Specifically, we use a combination of a logarithmic transformation and the half-vectorization operator to define an invertible transform $\hat{\gamma} = \vech\left(\ln(\hat{\Gamma})\right)$. Note that the half-vectorization operator is related to the normal vectorization operator by $\vech(X) = L_{n} \vecop(X)$ and $\vecop(X) = D_{n} \vech(X)$ where $D_n$ and $L_n$ are the duplication and elimination matrices respectively. The reason for this transformation is so that the transformed variable $\hat{\gamma}$ belongs to an unrestricted, unbounded set ($\hat{\gamma}\in\mathbb{R}^{n(n+1)/2}$) so that standard stability theorems, like those in \cite[Ch~4]{ref:khalil-2002}, can be used. To describe the dynamics of $\hat{\gamma}$, one needs to find the eigenvalue decomposition $\hat{\Gamma} = \Sigma_{\hat{\Gamma}} \Lambda_{\hat{\Gamma}} \Sigma_{\hat{\Gamma}}^T$, where $\Lambda_{\hat{\Gamma}} = \text{diag}\left(\lambda_1\left({\hat{\Gamma}}\right),\dots,\lambda_n\left({\hat{\Gamma}}\right)\right)$, to use the Dalecki{\v i}-Kre{\v i}n Theorem \cite[Th~2.10]{ref:norferini-2017}. Application of this theorem will result in
\begin{equation}
	\label{eq:nesc:model-free:alt-inverse-hessian}
	\frac{d}{dt} \hat{\gamma} = \vech\left(\Sigma_{\hat{\Gamma}}\left[C\left({\hat{\Gamma}}\right) \odot \left(\Sigma_{\hat{\Gamma}}^T\ \frac{d}{dt}\hat{\Gamma}\ \Sigma_{\hat{\Gamma}}\right)\right]\Sigma_{\hat{\Gamma}}^T\right)
\end{equation}
where
\begin{equation}
	C_{ij}\left({\hat{\Gamma}}\right) = \left\lbrace\begin{matrix}
		\frac{\ln\left(\lambda_i\left(\hat{\Gamma}\right)\right) - \ln\left(\lambda_j\left(\hat{\Gamma}\right)\right)}{\lambda_i\left(\hat{\Gamma}\right) - \lambda_j\left(\hat{\Gamma}\right)}&\hspace*{-8pt} \text{ if } \lambda_i\left(\hat{\Gamma}\right) \neq \lambda_j\left(\hat{\Gamma}\right) \\
		\frac{1}{\lambda_i\left(\hat{\Gamma}\right)} & \text{ otherwise }
	\end{matrix}\right.
\end{equation}
and $\odot$ is the Hadamard product. The average system and the model-based have the transformed states $\hat{\bar{\gamma}} = \vech\left(\ln\left(\hat{\bar{\Gamma}}\right)\right)$ and $\varpi = \vech\left(\ln\left(\Pi\right)\right)$ and the corresponding differential equations 
\begin{align}
    \label{eq:nesc:average-system:alt-inverse-hessian}
    \frac{d}{dt} \hat{\bar{\gamma}} &{}={} \vech\left(\Sigma_{\hat{\bar{\Gamma}}}\left[C\left({\hat{\bar{\Gamma}}}\right) \odot \left(\Sigma_{\hat{\bar{\Gamma}}}^T\ \frac{d}{dt}\hat{\bar{\Gamma}}\ \Sigma_{\hat{\bar{\Gamma}}}\right)\right]\Sigma_{\hat{\bar{\Gamma}}}^T\right) \\
    \label{eq:nesc:model-based:alt-inverse-hessian}
    \frac{d}{dt} \varpi &{}={} \vech\left(\Sigma_\Pi \left[C(\Pi) \odot \left(\Sigma_\Pi^T\ \frac{d}{dt}\Pi\ \Sigma_\Pi\right)\right]\Sigma_\Pi^T\right)
\end{align}
where $\Sigma_{\hat{\bar{\Gamma}}}$ and $\Sigma_{\Pi}$ are from the eigenvalue decompositions $\hat{\bar{\Gamma}} = \Sigma_{\hat{\bar{\Gamma}}} \Lambda_{\hat{\bar{\Gamma}}} \Sigma_{\hat{\bar{\Gamma}}}^T$ and $\Pi = \Sigma_{\Pi} \Lambda_{\Pi} \Sigma_{\Pi}^T$ respectively.

\subsection{Determining Model-Based ESCs from Model-Free ESCs without Calculating the Average Systems}
\label{sec:extremum-seekers:mb-from-mf}

Previous literature, such as \cite{ref:nesic-2010} and \cite{ref:nesic-2013}, started with model-based ESCs, often inspired by discrete-time optimizers and their continuous-time implementations, to construct the model-free ESCs. This leads to the practical stability of the model-free ESCs but often introduces a third small parameter which controls the speed of boundary layer dynamics. This 1) complicates the analysis and 2) adds another time scale which we find undesirable for the reasons discussed in Section~\ref{sec:intro}. Rather than using model-based ESCs as the inspiration for model-free ESCs and having three small parameters, we explore \emph{starting} with the model-free ESCs for generating associated model-based ESCs and only having, \emph{at most}, two small parameters. 

We leverage the fact that the average estimates in \eqref{eq:gesc:average-gradient-estimate} and \eqref{eq:nesc:average-hessian-estimate} converge to the true derivatives $\nabla J$ and $\nabla^2 J$ as $a\to 0$, provided $J$ is sufficiently differentiable and $s$, $m$, and $N$ are appropriately chosen by the system designer. Thus, the model-based GESC and NESC systems in \eqref{eq:gesc:model-based:parameter-update}, \eqref{eq:nesc:model-based:parameter-update-ode}, and \eqref{eq:nesc:model-based:inverse-hessian-ode} can be quickly obtained from the model-free GESC and NESC systems in \eqref{eq:gesc:model-free:parameter-update}, \eqref{eq:nesc:model-free:parameter-update-ode}, and \eqref{eq:nesc:model-free:inverse-hessian-ode} by replacing $m$ and $N$ with the operators $\nabla$ and $\nabla^2$, and removing the perturbations. This means that multiplying $m$ and $N$ by the perturbed map of $J$ is analogous to differentiating $J$ provided $J$ is sufficiently differentiable. The analogy between multiplying $m$ and $N$ and differentiation becomes closer as $a$ vanishes. However, since $J$ is unknown \emph{a priori}, this analogy to differentiation must be true for any function that is sufficiently differentiable in $\theta$, not just the $J$ in question. This means that we can regularize the vector field of the model-free ESC by taking nonlinear transforms of $J$, yet still know the associated model-based ESC. The following two propositions formalize this insight.
\begin{prop}
    \label{prop:scalar-compositions:gradient}
    Suppose $s$ and $m$ are almost periodic functions chosen such that \eqref{eq:gesc:average-gradient-estimate} is true for any continuously differentiable cost function $J:\mathbb{R}^n \to \mathcal{I} \subseteq \mathbb{R}$. If $v:\mathcal{I}\to\mathbb{R}$ is continuously differentiable, then the function composition $v\circ J$ satisfies
    \begin{equation}
        \label{eq:prop:scalar-compositions:gradient}
        \lim_{T\to\infty} \frac{1}{T} \int_0^T m(\tau, a)\cdot (v\circ J)\left(\hat{\bar{\theta}} + a s(\tau)\right) d\tau =
        \nabla \left(v \circ J\right) \left(\hat{\bar{\theta}}\right) + R_{v,m}\left(\hat{\bar{\theta}}, a\right)
    \end{equation}
    with $R_{v,m}\left(\cdot, a\right) \to 0$ compactly as $a\to 0$. 
\end{prop}
\begin{pf}
    The composition $v\circ J :\mathbb{R}^n \to \mathbb{R}$ has a gradient
    \begin{equation}
        \nabla (v\circ J)(\theta) = (v'\circ J)(\theta) \nabla J(\theta)
    \end{equation}
    which is well defined $\forall\theta \in \mathbb{R}^n$. Hence, $v \circ J$ has the first-order Taylor series expansion
    \begin{equation}
    (v\circ J)(\hat{\bar{\theta}} + \theta) = (v\circ J)(\hat{\bar{\theta}}) + (v'\circ J)(\hat{\bar{\theta}})\ \left( \nabla J(\hat{\bar{\theta}}) \cdot \theta\right) + R_{v,1}\left(\theta\right)
\end{equation}
where $R_{v,1}$ is the remainder of the first-order approximation at $\hat{\bar{\theta}}$. The rest of the proof follows the same averaging steps used to establish \eqref{eq:gesc:average-gradient-estimate}. These steps can be found in various works such as \cite{ref:krstic-2000,ref:ghaffari-2012,ref:kavanaugh-2025}. The convergence rate result from \cite{ref:kavanaugh-2025} and $v'\circ J$ being uniformly continuous on compact sets proves that $R_{v,m}(\cdot,a)\to 0$ compactly.
\end{pf}
\begin{prop}
    \label{prop:scalar-compositions:hessian}
    Suppose $s$ and $N$ are almost periodic functions chosen such that \eqref{eq:nesc:average-hessian-estimate} is true for any $\mathcal{C}^2$ cost function $J:\mathbb{R}^n \to \mathcal{I} \subseteq \mathbb{R}$. If $v:\mathcal{I}\to\mathbb{R}$ is $\mathcal{C}^2$, then the function composition $v\circ J$ satisfies
    \begin{equation}
        \label{eq:prop:scalar-compositions:hessian}
        \lim_{T\to\infty} \frac{1}{T} \int_0^T N(\tau, a) \cdot (v\circ J)(\hat{\bar{\theta}} + a s(\tau)) d\tau =
        \nabla^2 \left(v \circ J\right) (\hat{\bar{\theta}}) + R_{v,N}(\hat{\bar{\theta}}, a)
    \end{equation}
    with $R_{v,N}(\cdot, a) \to 0$ compactly as $a\to 0$.
\end{prop}
\begin{pf}
    The proof is the same as for Proposition~\ref{prop:scalar-compositions:gradient}, except we now use the second-order Taylor series expansion
    \begin{equation}
        v\circ J\left(\hat{\bar{\theta}} + \theta\right) = (v\circ J)\left(\hat{\bar{\theta}}\right) + \nabla (v\circ J)\left(\hat{\bar{\theta}}\right) \cdot \theta + \frac{1}{2} \theta^\top \nabla^2 (v\circ J)\left(\hat{\bar{\theta}}\right) \theta + R_{v,2}(\theta)
    \end{equation}
    where the Hessian $\nabla^2 (v\circ J)$ is defined as
    \begin{equation}
        \label{eq:scalar-compositions:hessian:chain-rule}
        \nabla^2 (v\circ J) \left(\hat{\bar{\theta}}\right) = (v' \circ J)\left(\hat{\bar{\theta}}\right)\ \nabla^2 J\left(\hat{\bar{\theta}}\right) + (v'' \circ J)\left(\hat{\bar{\theta}}\right)\ \nabla J\left(\hat{\bar{\theta}}\right) \nabla J\left(\hat{\bar{\theta}}\right)^{\top}
    \end{equation}
    by the chain rule and $R_{v,2}$ is the second-order remainder of the approximation of $v\circ J$ at $\hat{\bar{\theta}}$.
\end{pf}
  	
To give an example of the useful applications of Proposition~\ref{prop:scalar-compositions:gradient}, consider the dynamics of a model-free ESC given by
\begin{equation}
    \label{eq:ex-sys-1}
    \frac{d}{dt} \hat{\theta} = -k\ m(\omega t, a)\ \tanh \left(J\left(\hat{\theta} + a s(\omega t)\right)\right)
\end{equation}
with $\tanh\circ J$ being the function composition of interest. The reason to consider this model-free ESC is that the bounded output of the hyperbolic tangent bounds the speed of $\hat{\theta}$ by 
\begin{equation}
    \label{eq:ex-sys-1:speed-bound}
    \left\Vert \frac{d}{dt} \hat{\theta} \right\Vert \leq \frac{2 k}{a}\left(\sum_{i=1}^{n} \frac{1}{r_i^2}\right)^{1/2}
\end{equation}
which is independent of $\hat{\theta}$. This means that the example system does not have any finite-escape time trajectories for any initial $\hat{\theta}$, regardless of the chosen $\omega$ and $a$. Finite-time escape trajectories cannot be ruled out for traditional model-free GESC or NESC, as the magnitude of $\frac{d}{dt}\hat{\theta}$ depends on $\hat{\theta}$ through $J$. In fact, the example system offers a distinct advantage over classic GESC or NESC since the existence of finite-time escape trajectories leads to difficulties in tuning the small parameters of these classic ESCs. However, previous literature involving derivative demodulation, such as \cite{ref:tan-2007, ref:nesic-2010, ref:kavanaugh-2025}, only mentions the derivatives of $J$; therefore, deriving the model-based system of \eqref{eq:ex-sys-1} may appear impossible without knowledge of the differentiation equivalence. Fortunately, by knowing the differentiation equivalence, we can quickly give the associated model-based ESC as
\begin{equation}
    \label{eq:ex-sys-1:model-based}
    \frac{d}{dt} \vartheta = \frac{-k \nabla J(\vartheta)}{\cosh(J(\vartheta))^2}
\end{equation}
from Proposition~\ref{prop:scalar-compositions:gradient}. In theory, we still formally give the average system, as it is an important system in the sequence of system approximations and is used to determine the practical stability of the original system. However, in practice, one can often skip the tedious calculations of the average system and the limit $a\to 0$, as the intermediate steps yield a conclusion identical to the one reached via differentiation. 

Proposition~\ref{prop:scalar-compositions:hessian} also provides interesting modifications to the inverse Riccati filter given in \eqref{eq:nesc:model-free:inverse-hessian-ode} by replacing $J$ with $v\circ J$. If we choose $v$ in Proposition~\ref{prop:scalar-compositions:hessian} such that the first and second derivatives are strictly positive over the image of $J$, then the filter input remains a strictly positive definite matrix. However, this introduces a nonlinear weighting between the Hessian and gradient outer product of $J$. Inverses of low-passed gradient outer products are generally found in adaptive gradient methods, specifically AdaGrad \cite{ref:duchi-2011}, and are used to increase the sensitivity to directions which have not been previously explored along the trajectories. This means that $v\circ J$ in \eqref{eq:nesc:model-free:inverse-hessian-ode} results in a model-based ESC which attempts to move towards the center of a local quadratic approximation while still being sensitive to unexplored directions. This behavior can be tuned for $J$, particularly if the system designer knows what $J$ is relative to some unknown optimum $\theta^*$, as is often the case in grey-box optimization.

In addition to taking nonlinear transforms of $J$, one can also introduce auxiliary variables $\hat{\phi}\in\mathbb{R}^{n'}$ to the mappings of $J$ for the purpose of generating other model-free systems. The formality of this is stated in the following propositions.
\begin{prop}
    \rev{\label{prop:auxiliary-variables}}
    Suppose $J:\mathbb{R}^n \to \mathcal{I}\subseteq \mathbb{R}$ is continuously differentiable and the function  $w:\mathbb{R}^{n'}\times\mathcal{I}\to\mathbb{R}$ is a continuously differentiable function. If $s$ and $m$ are almost periodic functions chosen such that \eqref{eq:gesc:average-gradient-estimate} is true, then 
    \begin{equation}
        \lim_{T\to\infty} \frac{1}{T} \int_0^T m(\tau, a) w\left(\hat{\bar{\phi}}, J\left(\hat{\bar{\theta}} + a s(\tau)\right)\right) d\tau = \left.\frac{\partial w\left(\hat{\bar{\phi}}, y\right)}{\partial y}\right\vert_{y = J\left(\hat{\bar{\theta}}\right)} \nabla J\left(\hat{\bar{\theta}}\right) + R_{w,m}\left(\hat{\bar{\theta}}, \hat{\bar{\phi}}, a\right)
    \end{equation}
    with $R_{w,m}\left(\cdot, \cdot, a\right)\to 0$ compactly as $a\to 0$. 
\end{prop}
\begin{pf}
    The proof is the same as for the previous two propositions. One needs to first calculate the first-order Taylor series centered at $\left(\hat{\bar{\phi}}, \hat{\bar{\theta}}\right)$,
    \begin{multline}
         w\left(\hat{\bar{\phi}} + \phi, J\left(\hat{\bar{\theta}} + \theta\right)\right) = w\left(\hat{\bar{\phi}}, J\left(\hat{\bar{\theta}}\right)\right) + \sum_{i=1}^{n'} \frac{\partial w\left(\hat{\bar{\phi}}, y\right)}{\partial\hat{\bar{\phi}}_i}\Bigg\vert_{y=J\left(\hat{\bar{\theta}}\right)} \phi_i \\ + \frac{\partial w\left(\hat{\bar{\phi}}, y\right)}{\partial y} \Bigg\vert_{y=J\left(\hat{\bar{\theta}}\right)} \sum_{j=1}^{n} \frac{\partial J\left(\hat{\bar{\theta}}\right)}{\partial\hat{\bar{\theta}}_j} \theta_j + R_{w,1}\left(\phi, \theta\right)
    \end{multline}
    where $R_{w,1}$ is the remainder term. The rest of the proof then proceeds with the usual averaging steps. Since $\hat{\bar{\phi}}$ is considered ``frozen'' in the averaging process and the dither signal does not affect $\hat{\phi}$, the terms associated with a change in $\hat{\bar{\phi}}$ are zero. Additionally, because $m$ is a zero-mean vector, the partial derivative terms of $w$ with respect to $\hat{\bar{\phi}}$ do not contribute to the average system. This completes the proof.
\end{pf}
While it may seem odd to include auxiliary variables to parameterize the mappings, the simplest example of such a function $w$ has already been introduced in this work. One can consider \eqref{eq:nesc:model-free:parameter-update-ode} from NESC to be an example by setting $\hat{\phi} = \vech(\hat{\Gamma})$ and $w(\hat{\phi}, y) = -k \vech^{-1}(\hat{\phi}) y$. However, for a more complicated, novel example, assume that we know $J(\theta) \geq 0$ for all $\theta\in\mathbb{R}$ and then consider
\begin{equation}
    \label{eq:determining-mb:w:example-2:mf-esc}
    \frac{d}{dt} \hat{\theta} = -k\ m(\omega t, a) e^{-\hat{\phi}}\sqrt{1 + 2 e^{\hat{\phi}}J(\hat{\theta} + a s(\omega t))}
\end{equation}
as the parameter portion of a model-free ESC. The associated model-based system is
\begin{equation}
    \label{eq:determining-mb:w:example-2:mb-esc}
    \frac{d}{dt} \vartheta = -k \frac{J'(\vartheta)}{\sqrt{1 + 2 e^{\varphi}J(\vartheta)}}
\end{equation}
where $\varphi$ represents the auxiliary variables in the model-based ESC. Again, it was not necessary for us to calculate both the average system and then the limit $a\to 0$ to determine the model-based ESC system. Instead, to determine \eqref{eq:determining-mb:w:example-2:mb-esc} from \eqref{eq:determining-mb:w:example-2:mf-esc}, we only need to follow three simple steps: 1) replace the signal $m$ with the operator $\nabla$; 2) set $a$ to $0$ to remove the additive perturbation; and 3) replace $\hat{\theta}$ and $\hat{\phi}$ with $\vartheta$ and $\varphi$, respectively.

\subsection{The ESC Abstraction}
    \label{sec:extremum-seekers:abstraction}

We have given an incomplete list of known ESC systems in Sections \ref{sec:extremum-seekers:gesc} and \ref{sec:extremum-seekers:nesc} but we wish to show how these systems can be abstracted so that a general stability theorem can be applied to both of them and other model-free continuous-time ESC algorithms, such as those we gave in Section~\ref{sec:extremum-seekers:mb-from-mf}. We write the model-free ESC as
\begin{equation}
	\label{eq:perturbation-esc}
	\frac{d}{dt} \hat{x} = f(\omega t, \hat{x}, a)
\end{equation}
where $\hat{x}\in\mathbb{R}^m$ and $m\geq n$ is the augmented space with both the parameter estimates $\hat{\theta}$ as well as any auxiliary states, e.g.,\ $\hat{\gamma}$ in the NESC algorithm. If the system designer has constructed $f$ such that $f$ is Lipschitz in $\hat{x}$ for fixed $a$ and $f(\cdot, \hat{x}, a)$ is a periodic (or almost periodic) function, then the model-free ESC has an average system 
\begin{equation}
	\label{eq:average-system}
	\frac{d}{dt} \hat{\bar{x}} = \bar{f}(\hat{\bar{x}}, a)
    : = \lim_{T\to\infty} \frac{1}{T} \int_{0}^T f(\tau, \hat{\bar{x}}, a) d\tau\,.
\end{equation}
by \cite[Lem~4.6.4]{ref:sanders-2007}. Note that while we choose to integrate over the interval $[0,T]$, integrating over $[t_0, t_0 + T]$ for any $t_0\in\mathbb{R}$ results in the same $\bar{f}(\hat{\bar{x}}, a)$ and the limit converges uniformly in $t_0$ \cite[Cor~3.2]{ref:fink-1975}.

From the average system, we define the model-based ESC as
\begin{equation}
	\label{eq:perturbation-free-esc}
	\frac{d}{dt} z = h(z)
    := \lim_{a\rightarrow 0^+} \bar{f}(z, a) 
\end{equation}
where the adjective ``model-based'' is used because $h$ evolves the augmented state vector $z$ based on the true derivative values of $J$ with respect to $z$, such as the gradient $\nabla J(z)$ or the Hessian $\nabla^2 J (z)$. The term is used here for the perturbation-free, derivative-driven limit of the design at hand, and not for a pre-existing optimizer: as Section~\ref{sec:extremum-seekers:mb-from-mf} illustrates, $h$ need not be a continuous-time analogue of any known discrete-time method.

Before stating the main result we record the hypotheses under which the three systems \eqref{eq:perturbation-esc}, \eqref{eq:average-system}, and \eqref{eq:perturbation-free-esc} are related. They are of two kinds: regularity and uniformity of the average in the fast time $\tau$, which is what general averaging requires, and locally uniform convergence of the average family as $a\to 0^+$, which is what the demodulation structure of a model-free ESC supplies.
\begin{asmp}
    \label{asmp:averaging}
    There is $\bar{a} > 0$ such that the following hold.
    \begin{enumerate}
        \item[(A1)] \emph{Regularity:} $f$ is continuous on $\mathbb{R}\times\mathbb{R}^m\times(0,\bar{a}]$ and, for each fixed $a\in(0,\bar{a}]$ and each compact $\mathcal{D}\subset\mathbb{R}^m$, $f(\tau,\cdot,a)$ is Lipschitz on $\mathcal{D}$ uniformly in $\tau$ and $\sup_{\tau\in\mathbb{R},\ \hat{x}\in\mathcal{D}}\Vert f(\tau,\hat{x},a)\Vert < \infty$.
        \item[(A2)] \emph{Strong average:} For each fixed $a\in(0,\bar{a}]$, the limit in \eqref{eq:average-system} exists, and the convergence is uniform in the initial fast phase and locally uniform in the state: for each compact $\mathcal{D}\subset\mathbb{R}^m$ there exist $\kappa_{\mathcal{D}}(a) > 0$ and $\sigma\in\mathcal{L}$ such that
        \begin{equation}
            \label{eq:strong-average}
            \sup_{\tau_0\in\mathbb{R}}\ \sup_{\hat{x}\in\mathcal{D}}\left\Vert \frac{1}{T}\int_{\tau_0}^{\tau_0 + T}\left[f(\tau,\hat{x},a) - \bar{f}(\hat{x},a)\right] d\tau \right\Vert \leq \kappa_{\mathcal{D}}(a)\,\sigma(T)\,,\quad \forall\, T > 0\,.
        \end{equation}
        \item[(A3)] \emph{Compact convergence:} $\bar{f}(\cdot,a)\to h(\cdot)$ uniformly on compact subsets of $\mathbb{R}^m$ as $a\to 0^+$.
        \item[(A4)] \emph{Stability of the limit:} $h(x^*) = 0$, every maximal solution of \eqref{eq:perturbation-free-esc} is forward complete, and \eqref{eq:perturbation-free-esc} is GAS to $x^*$ in the sense that the stability and attractivity estimates hold for \emph{every} solution of the initial value problem, solutions being possibly nonunique.
    \end{enumerate}
\end{asmp}
Conditions (A1) and (A2) are the Krylov--Bogoliubov--Mitropolsky requirements under which the general averaging theory of \cite[Ch~4]{ref:sanders-2007} applies, stated so that the constants are uniform on the compact sets and initial phases used in the induction of Section~\ref{sec:main-results:model-free}. They are met by the ESC constructions of this paper: for fixed $a$ the vector fields are continuous and locally Lipschitz in $\hat{x}$ whenever $J\in\mathcal{C}^1$, and they are almost periodic in $\tau$ uniformly for $\hat{x}$ in compact sets, since $J$ is uniformly continuous on compact sets and $s$ is Bohr almost periodic, which yields \eqref{eq:strong-average} through \cite[Th~1.9, Th~3.1, Cor~3.2]{ref:fink-1975}. Condition (A3) is the compact convergence established for the GESC after \eqref{eq:gesc:model-based:parameter-update} and, through Propositions~\ref{prop:scalar-compositions:gradient}--\ref{prop:auxiliary-variables}, for the designs derived in Section~\ref{sec:extremum-seekers:mb-from-mf}. Condition (A4) is the form of GAS required by the converse Lyapunov theorem invoked in Section~\ref{sec:main-results:average-system}, in which uniqueness of solutions of \eqref{eq:perturbation-free-esc} is not assumed \cite{ref:kurzweil-1956,ref:kellett-2015}.

We now give the main result of this work where stability definitions in the theorem are given in Appendix \ref{sec:practical-stability-definitions} for reference. 
	
\begin{thm}
	\label{thm:model-based-implies-model-free}
	Consider the model-free, perturbation-based ESC in \eqref{eq:perturbation-esc} and its associated average system \eqref{eq:average-system} and model-based, perturbation-free ESC \eqref{eq:perturbation-free-esc}. Let Assumption~\ref{asmp:averaging} hold with $x^*$ the augmented state vector which corresponds to $\theta^*$ and, without loss of generality, taken to be the origin. Then
	\begin{enumerate}
        \label{thm:model-based-implies-model-free:average-system}
		\item the average system defined in \eqref{eq:average-system} is semiglobally practically asymptotically stable (sGPAS) to the origin with respect to the small parameter $a$ and
        \label{thm:model-based-implies-model-free:model-free}
		\item the model-free ESC system is sGPUAS to the origin with respect to the small parameter vector $(a,\omega^{-1})$.
	\end{enumerate}
    In particular, for any given pair of initial and ultimate radii, the trajectories of \eqref{eq:perturbation-esc} can be made by choice of the small parameter vector $(a,\omega^{-1})$ to remain in a compact set for all $t\geq t_0$; thus, \eqref{eq:perturbation-esc} is forward complete from the initial conditions considered.
\end{thm}

\begin{rem}
	\em
    \label{rem:sharpness}
    The uniformity in the initial fast phase demanded by \eqref{eq:strong-average} cannot be weakened to existence of the limit in \eqref{eq:average-system} for a single $\tau_0$.\footnote{The authors thank an anonymous reviewer of an earlier version of the manuscript for a comment that gives rise to this remark.} Take $m=1$, with $a$ absent from the dynamics, and
    \begin{equation}
        f(\tau,\hat{x},a) = -\hat{x} + u(\tau)\,,\qquad
        u(\tau) = \begin{cases}
            1 & \tau\in[k^3,\ k^3 + k) \\
            0 & \tau\in[k^3 + k,\ (k+1)^3)
        \end{cases}
    \end{equation}
    for $k = 0,1,2,\ldots$. The fraction of $[0,N^3]$ on which $u(\tau) = 1$ is of order $\mathcal{O}(N^{-1})$; thus, the limit in \eqref{eq:average-system} exists with $\tau_0 = 0$ and gives $\bar{f}(\hat{\bar{x}},a) = h(\hat{\bar{x}}) = -\hat{\bar{x}}$, which is GAS. The conclusions of Theorem~\ref{thm:model-based-implies-model-free} nevertheless fail: on the $k$th block the state is driven toward $1$ for a fast-time duration $k$ that grows without bound, so no choice of $\omega$ confines the trajectories of \eqref{eq:perturbation-esc} to a small neighborhood of the origin. What fails is \eqref{eq:strong-average}, whose left-hand side equals $1$ for every $k$ upon taking $\tau_0 = k^3$ and $T = k$. Almost periodicity of the dither rules out such behavior, and this is the role it plays in Assumption~\ref{asmp:averaging}.
\end{rem}

The primary difference between the above result and those in previous literature is in the assumptions. Previous literature, such as \cite{ref:nesic-2013}, has generally assumed that the model-based system is some continuous-time analogue of a well-known, discrete-time optimization algorithm where derivatives are treated as inputs to the model-free system. As such, it was then necessary to verify that the continuous-time optimization algorithm had certain robustness properties with regard to the magnitude of additive disturbances to the derivative inputs (cf.~\cite[Prop~2]{ref:nesic-2010}). In contrast, this work only makes assumptions about the model-free system and the cost function, eliminating the need to verify the aforementioned property for the model-based system. Rather, the practical stability results are a natural consequence of how model-free systems are commonly constructed. The use of quickly perturbed $J$ multiplied by $m$ or $N$ in the model-free systems naturally gives rise to average and model-based systems with the average system converging compactly to the model-based system. We demonstrated this in Section~\ref{sec:extremum-seekers:mb-from-mf} by showing interesting model-free ESCs that are not continuous-time analogues of well-known discrete-time optimization algorithms but nevertheless have average and model-based systems. Thus, the results in this work rely on simpler assumptions and can be used in a broader capacity.
					
\section{Motivating Example}
	
To illustrate the usefulness of Theorem \ref{thm:model-based-implies-model-free}, consider a GESC acting on the quartic, two-dimensional cost function
\begin{equation}
    \label{eq:motivating-example:cost-function}
	J(\theta) = \theta_1^4 + \left(\theta_1 + \theta_2\right)^4
\end{equation}
This $J$ is differentiable, has a unique minimum at the origin ($\theta^* = 0$), and a zero matrix Hessian at this minimum. Since $J$ is a polynomial of only quartic terms, there $\exists b_1,b_2 > 0$ to bound the minimum and maximum growth rates of $J$ by
\begin{equation}
	b_1 \left\Vert \theta \right\Vert_2^4 \leq J(\theta) \leq b_2 \left\Vert \theta \right\Vert_2^4\,.
\end{equation}
Clearly, this $J$ has the same properties of a Lyapunov candidate function. When considering the model-based GESC defined in \eqref{eq:gesc:model-based:parameter-update} acting on $J$, the Lie derivative of $J$ is
\begin{align}
	\label{eq:gesc-example:model-based:lyapunov-derivative}
	\frac{dJ}{dt}\left(\vartheta\right) &{}={} -k \left\Vert \nabla J \left(\vartheta \right) \right\Vert_2^2 
    {}={} -16k\left[\left(\vartheta_1^3 + (\vartheta_1 + \vartheta_2)^3\right)^2 + (\vartheta_1 + \vartheta_2)^6\right]
\end{align}
which is a negative definite function. Thus, the model-based GESC is GAS to the origin \cite[Th~4.9]{ref:khalil-2002} and, consequently, \emph{any} average system and model-free GESC that are associated with the model-based GESC are sGPAS and sGPUAS, respectively, by Theorem \ref{thm:model-based-implies-model-free}.

However, consider attempting to prove just the result that the average system is sGPAS \emph{without} Theorem \ref{thm:model-based-implies-model-free}. Indeed, the stability analysis of the GESC's average system yields a much more nuanced story. Consider the dither signal $s$ with the chosen parameters $\omega'_1 = 1$ and $\omega'_2 = 3$ for the relative dither rates. These dither signals result in the average gradient estimate of 
\begin{equation}
	\label{eq:gesc-example:average-system}
	\hat{\bar{g}}(\hat{\bar{\theta}}, a) = 4\begin{bmatrix}
		\hat{\bar{\theta}}_1^3 + \left(\hat{\bar{\theta}}_1 + \hat{\bar{\theta}}_2\right)^3 \\ \left(\hat{\bar{\theta}}_1 + \hat{\bar{\theta}}_2\right)^3
	\end{bmatrix} + a^2 A(r_1, r_2) \hat{\bar{\theta}}
\end{equation}
where 
\begin{equation}
	\label{eq:gesc-example:A}
	A(r_1, r_2) = \begin{bmatrix}
		6 r_1^2 - 3 r_1 r_2 + 6 r_2^2 & 3 r_1^2 - 3 r_1 r_2 + 6 r_2^2 \\
		\frac{-2 r_1^3}{r_2} + 6 r_1^2 + 3 r_2^2 & \frac{- r_1^{3}}{r_2} + 6 r_1^2 + 3 r_2^2
	\end{bmatrix}
\end{equation}
The cubic terms in \eqref{eq:gesc-example:average-system} are $\nabla J$ and the linear terms comprise $R_g$. 

\begin{figure}[t!]
	\centering
	\includegraphics[width=3in]{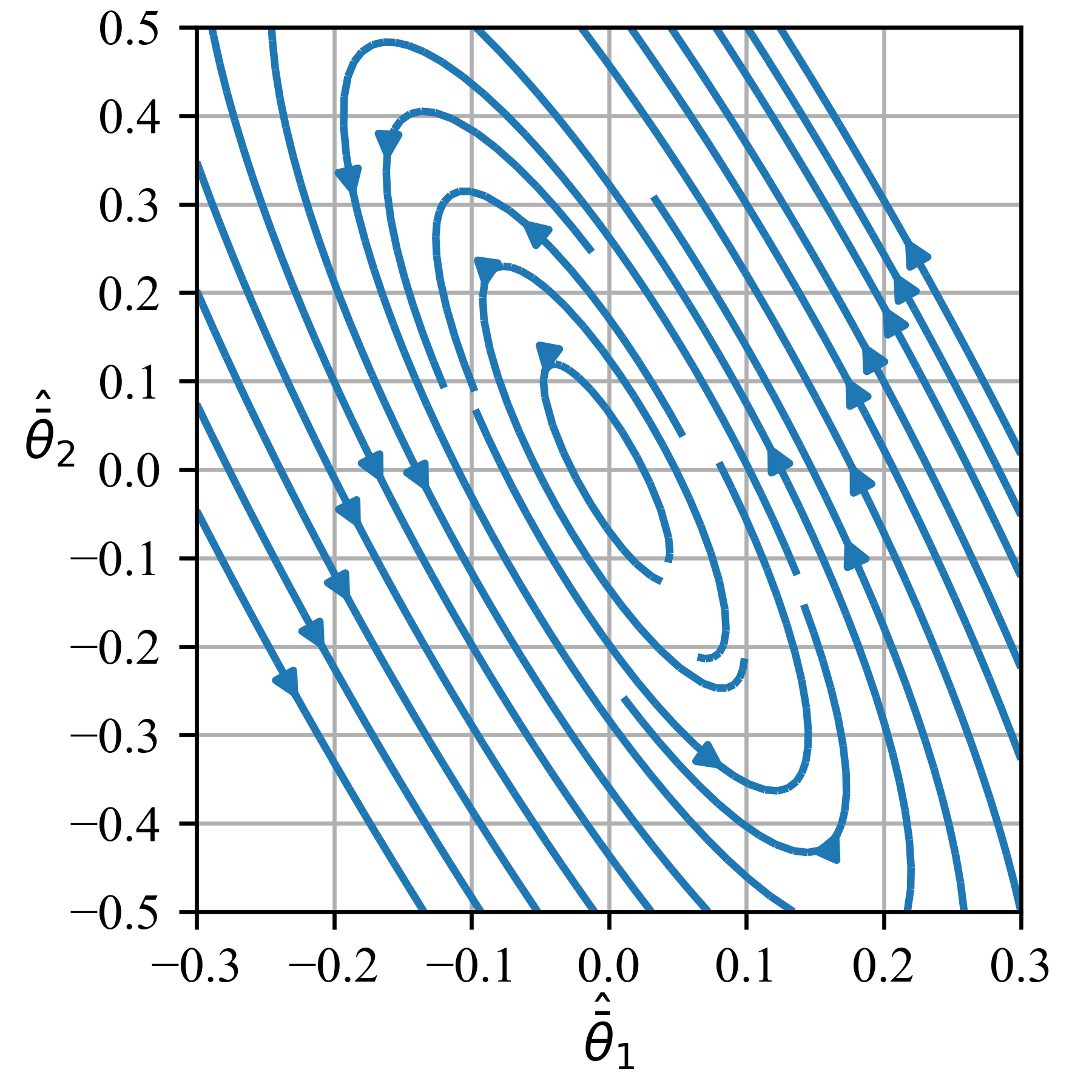}\hfill
	\includegraphics[width=3in]{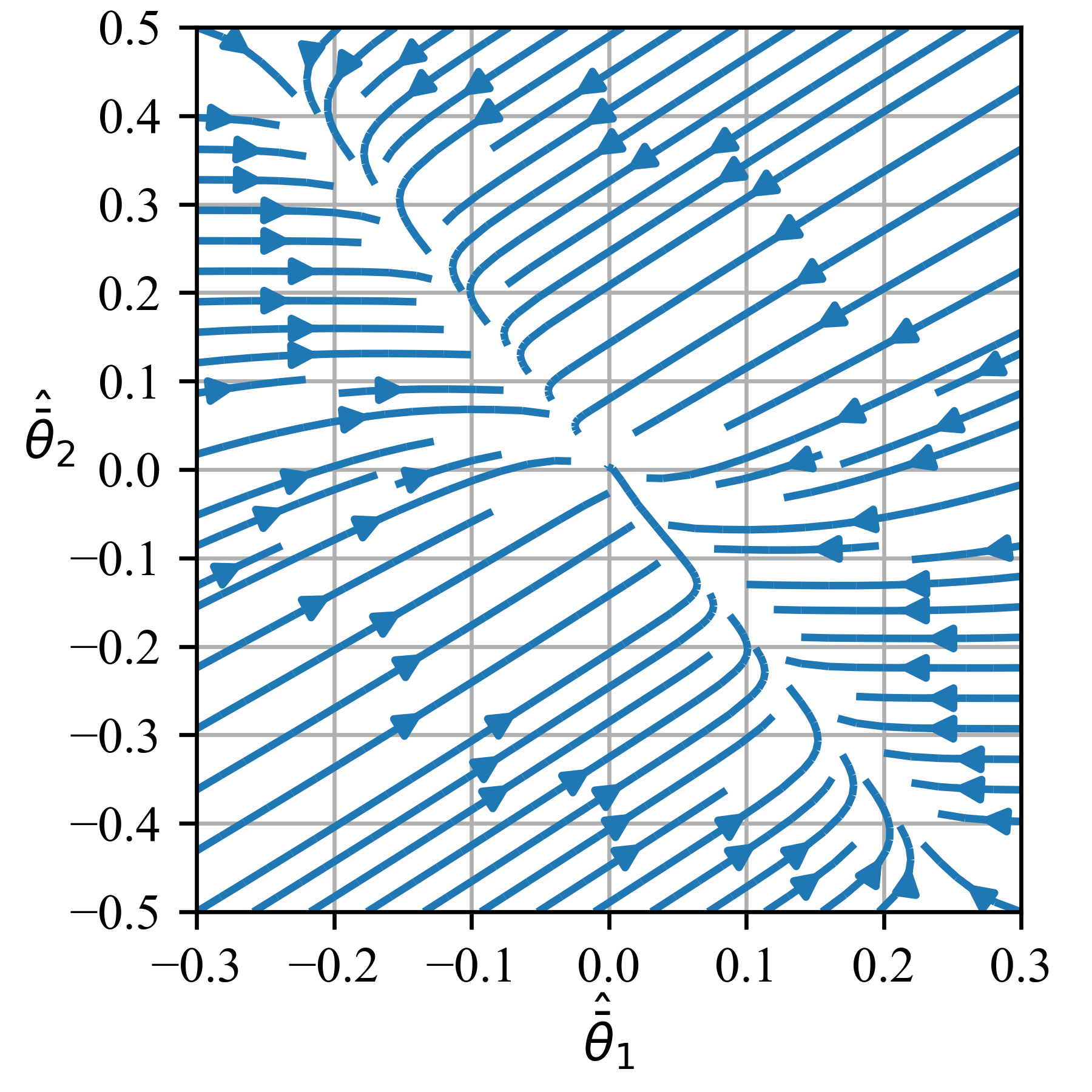}
	\caption{Stream plots of the average GESC system from the considered example. Left shows results with $a=100$ and right shows results with $a=0.01$. The origin is an unstable focus in both cases, and each carries a limit cycle of radius $0.2427\,a$; by \eqref{eq:motivating-example:self-similar} the two panels differ only by the magnification $10^4$ at which the same portrait is viewed. Over the fixed window shown, the cycle for $a=100$ lies far outside the frame while the cycle for $a=0.01$ is too small to resolve, which is what practical stability with respect to $a$ amounts to here.}
    \label{fig:example:trajectories}
\end{figure}

The origin is an equilibrium of the average system for any choice of $r_1$ and $r_2$ just as it was for the model-based ESC. However, the local stability of this equilibrium is based on the linear terms, which are associated only with $R_g$ as $\nabla^2 J(0)=0$. When we linearize the system, the differential equation is
\begin{equation}
	\frac{d\hat{\bar{\theta}}}{dt} = -k a^2 A(r_1, r_2) \hat{\bar{\theta}}
\end{equation}
For local stability, a necessary condition is that the eigenvalues of the Jacobian have strictly negative real parts. Consider the specific relative dither amplitudes of $r_1 = 12/\sqrt{145}$ and $r_2 = 1/\sqrt{145}$ where the specific linearized system is
\begin{equation}
	\frac{d\hat{\bar{\theta}}}{dt} = -\frac{k a^2}{145} \begin{bmatrix}
		834 & 402 \\ -2589 & -861
	\end{bmatrix} \hat{\bar{\theta}}
\end{equation}
The two complex eigenvalues of this system
\begin{equation}
	\lambda_{1,2} = \frac{9 k a^2}{290} \left(3 \pm i \sqrt{15927}\right)
\end{equation}
have positive real parts, which implies that the average system is locally unstable and consequently not GAS. This seems to be an impasse; if the average system is not GAS, then why should the average system be sGPAS and the model-free ESC be sGPUAS? The reason why is that practical stability is a weaker notion of stability than traditional stability notions. Although the origin is unstable for the average system when $a > 0$, it is still sGPAS to the origin. To illustrate this, examine the trajectories in Fig.~\ref{fig:example:trajectories}. The origin is an unstable focus for \emph{every} $a>0$, and the trajectories approach a stable limit cycle whose diameter shrinks with $a$. In this example the dependence on $a$ is exact rather than merely qualitative. Since $\nabla J$ is homogeneous of degree three and $R_g = a^2 A(r_1,r_2)\hat{\bar{\theta}}$ is linear in the state, the substitution $\hat{\bar{\theta}} = a\xi$ turns the average system into
\begin{equation}
    \label{eq:motivating-example:self-similar}
    \frac{d}{dt}\xi = -k a^3 \left(4\begin{bmatrix} \xi_1^3 + (\xi_1 + \xi_2)^3 \\ (\xi_1 + \xi_2)^3\end{bmatrix} + A(r_1,r_2)\,\xi\right)\,,
\end{equation}
in which $a$ appears only as the time scale $ka^3$. The phase portrait in $\xi$ is therefore independent of $a$: the attractor is one fixed curve in $\xi$, and the limit cycle in $\hat{\bar{\theta}}$ is that curve scaled by $a$, so its radius is exactly proportional to the dither amplitude. Numerically, the cycle for these $(r_1,r_2)$ has radius $0.2427\,a$ to four digits over decades of $a$. The two panels of Fig.~\ref{fig:example:trajectories} are consequently the same portrait viewed at two magnifications rather than two qualitatively different behaviors, and sGPAS for this example follows from \eqref{eq:motivating-example:self-similar} alone. We give the Lyapunov argument anyway, since it is the one that generalizes: homogeneity is a property of this particular $J$, whereas the estimate below only uses the growth bounds on $J$. What we need to show is that the system designer should be able to choose a sufficiently small $a$ such that the system meets the required convergence criteria.
	
Consider the Lyapunov function $V$ defined as $V = \frac{1}{2} \left\Vert \hat{\bar{\theta}} \right\Vert_2^2$, then the Lie derivative is
	\begin{align}
		\label{eq:motivating-example:lyapunov-derivative:full}
		\frac{dV}{dt} &{}={} -k\left[\left( 8 \hat{\bar{\theta}}_1^4 + 16 \hat{\bar{\theta}}_1^3 \hat{\bar{\theta}}_2 + 24 \hat{\bar{\theta}}_1^2 \hat{\bar{\theta}}_2^2 + 16 \hat{\bar{\theta}}_1 \hat{\bar{\theta}}_2^3 + 4 \hat{\bar{\theta}}_2^4 \right) + \frac{a^2}{145}\left(834\hat{\bar{\theta}}_1^2 - 2187\hat{\bar{\theta}}_1 \hat{\bar{\theta}}_2 - 861\hat{\bar{\theta}}_2^2\right)\right] \\
		\label{eq:motivating-example:lyapunov-derivative:simplified}
		&{}={} -4k J(\hat{\bar{\theta}}) - \frac{ka^2}{145} \left(834\hat{\bar{\theta}}_1^2 - 2187\hat{\bar{\theta}}_1 \hat{\bar{\theta}}_2 - 861\hat{\bar{\theta}}_2^2\right)
	\end{align}
The $-4kJ(\hat{\bar{\theta}})$ component is contributed by $\nabla J$ whereas the terms associated with $a^2$ are contributed by $R_g$. We then bound the Lie derivative by
\begin{align}
	\frac{dV}{dt} &\leq -4k J(\hat{\bar{\theta}}) + 18 k a^2 \left\Vert \hat{\bar{\theta}} \right\Vert_2^2 
	\leq -4k b_1 \left\Vert \hat{\bar{\theta}} \right\Vert_2^4 + 18 k a^2 \left\Vert \hat{\bar{\theta}} \right\Vert_2^2
\end{align}
since $\hat{\bar{\theta}}^\top A \hat{\bar{\theta}} \geq -\lambda_1\left(A_s(r_1,r_2)\right)\Vert\hat{\bar{\theta}}\Vert_2^2$, where $A_s = \frac{1}{2}\left(A + A^\top\right)$ and $\lambda_1$ denotes the largest eigenvalue in magnitude with $\lambda_1(A_s) < 9.64$ for these relative amplitudes. From this upper bound, it can be shown that the condition
\begin{equation}	
	\left\Vert \hat{\bar{\theta}} \right\Vert_2 > \frac{3a}{\sqrt{b_1}} \implies \frac{dV}{dt} \left(\hat{\bar{\theta}}\right) < -2 k b_1 \left\Vert \hat{\bar{\theta}} \right\Vert_2^4
\end{equation} 
The choice of $V$ and the condition for $\frac{dV}{dt} < 0$ appears similar to conditions laid out in \cite[Th~4.19]{ref:khalil-2002} for proving Input-to-State Stability (ISS) of a system if $a$ is somehow considered some sort of disturbance to the system. Hence we can bound the trajectories of the average system by
\begin{equation}
	\left\Vert \hat{\bar{\theta}}(t) \right\Vert \leq \beta\left(\left\Vert \hat{\bar{\theta}}(t_0) \right\Vert, t - t_0\right) + \gamma(a),\ \forall t\in[t_0,\infty)
\end{equation}
where $\beta\in\mathcal{KL}$ and $\gamma(a) = 3 b_1^{-3/4} b_2^{1/4} a$. Note that the ``disturbance" $a$ is \emph{a priori} known to the system designer, who sets the absolute magnitude of the dither signal for the deployment of the model-free ESC. Thus for any $c_2> 0$, all $a\in(0, c_2 b_1^{3/4}/3 b_2^{1/4})$ result in trajectories of the average system satisfying
\begin{equation}
	\left\Vert \hat{\bar{\theta}}(t) \right\Vert \leq \beta\left(\left\Vert \hat{\bar{\theta}}(t_0) \right\Vert, t - t_0\right) + c_2,\ \forall t\in[t_0,\infty)
\end{equation}
which is the definition of global practical asymptotic stability (GPAS), and therefore sGPAS, given in Definition~\ref{def:sgpuas:alternative-definition} in Appendix \ref{sec:practical-stability-definitions}.

We can expand this insight by considering other dither signals $s$ with any other valid pair $(r_1,r_2)$, where neither $r_1$ nor $r_2$ is zero and $r_1^2 + r_2^2 = 1$. In this more general average system, the Lie derivative of $V$ can be upper bounded by
\begin{equation}
	\frac{d}{dt} V(\hat{\bar{\theta}}) \leq - 4 k b_1 \left\Vert \hat{\bar{\theta}} \right\Vert_2^4 + k a^2 \lambda_1 \left( A_s(r_1,r_2)\right) \left\Vert \hat{\bar{\theta}} \right\Vert_2^2
\end{equation}
Again, there is a condition
\begin{equation}
	\left\Vert \hat{\bar{\theta}} \right\Vert_2 > a\sqrt{\frac{\lambda_1(A_s(r_1,r_2))}{2b_1}} \implies \frac{dV}{dt} \leq -2k b_1 \left\Vert \hat{\bar{\theta}} \right\Vert_2^4
\end{equation}
which will prove the ISS of the average system with respect to the ``disturbance" $a$ \cite[Th~4.19]{ref:khalil-2002}. Hence the trajectories of the average system for this $(r_1,r_2)$ pair have the bound
\begin{equation}
	\label{eq:gesc-example:d1-d2-bound}
	\left\Vert \hat{\bar{\theta}}(t) \right\Vert \leq \beta\left( \left\Vert \hat{\bar{\theta}}(t_0) \right\Vert, t-t_0\right) + a \sqrt[4]{\frac{b_2}{b_1}} \sqrt{\frac{\lambda_1(A_s(r_1,r_2))}{2 b_1}}
\end{equation}
and as $a$ can be made arbitrarily small, the average system for this $(r_1,r_2)$ pair is GPAS with respect to the small parameter $a$.

\begin{figure}[t!]
    \centering
    \includegraphics[width=3.2in]{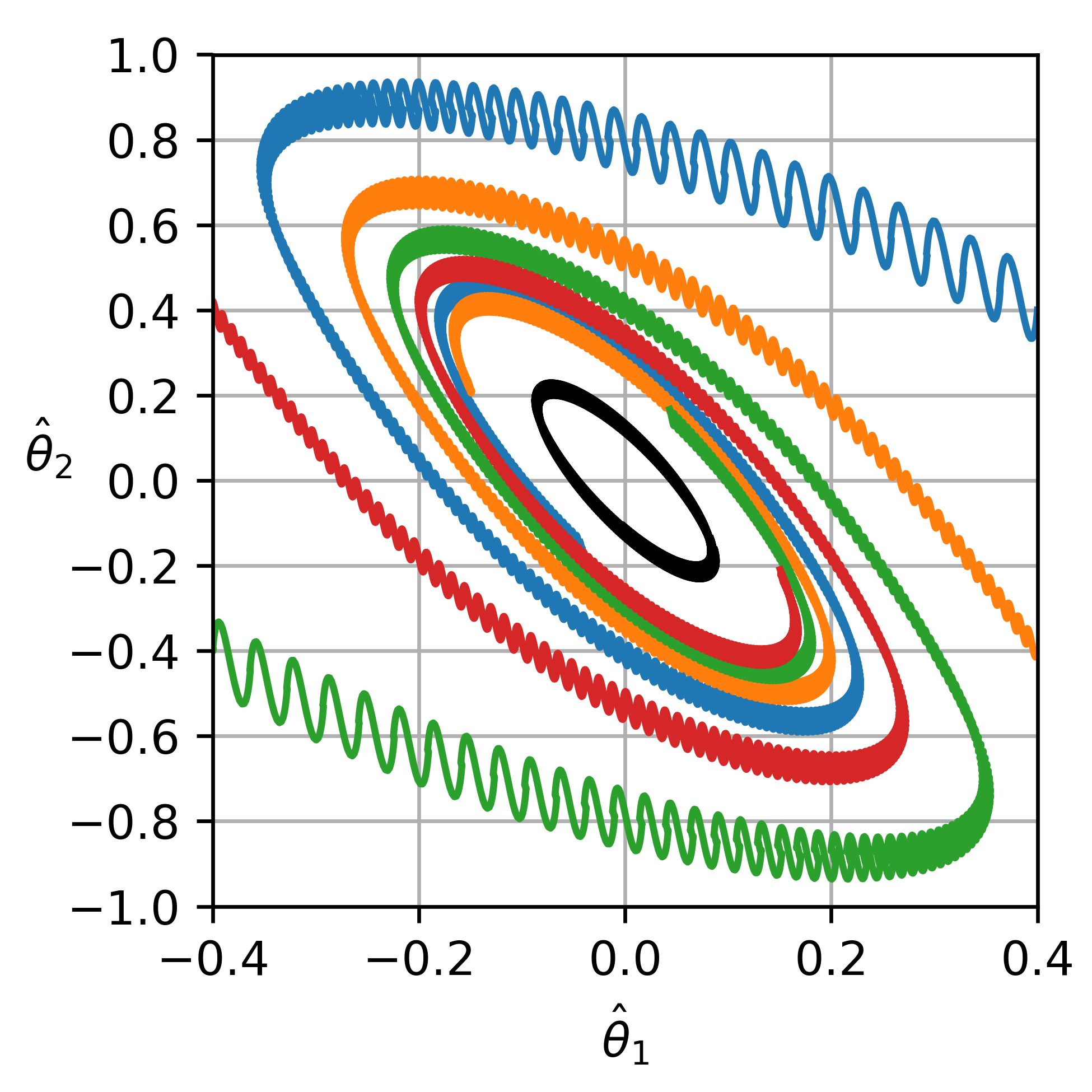}\hfill
    \includegraphics[width=3.2in]{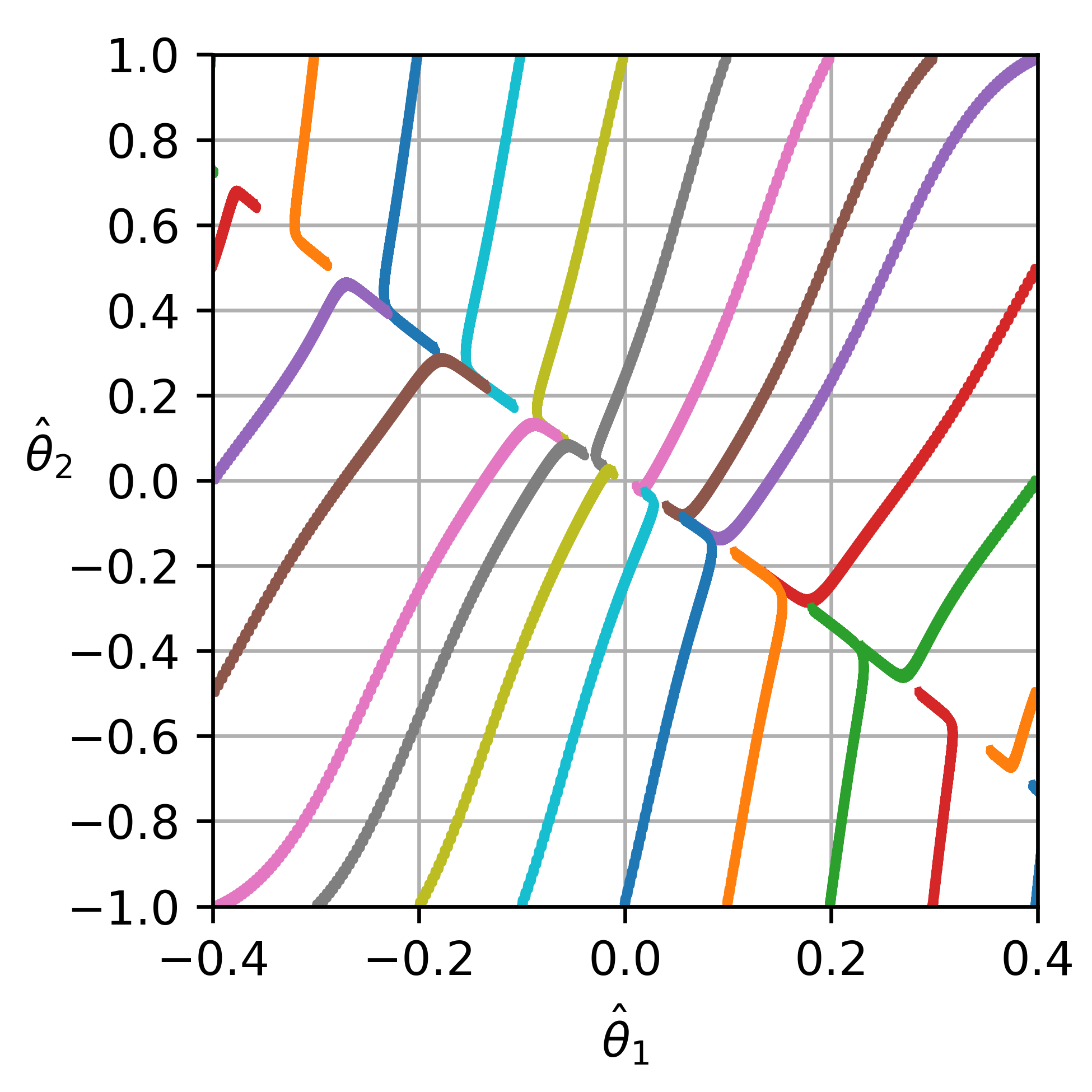}
    \caption{Simulations of the model-free ESCs with the cost function in \eqref{eq:motivating-example:cost-function} for 200 seconds. The left is the GESC in \eqref{eq:gesc:model-free:parameter-update} while the right is the novel ESC in \eqref{eq:ex-sys-1}. The simulation parameters used for both model-free ESCs were $k=0.01$, $\omega=10$, $\omega_1' = 1$, $\omega'_2 = 3$, $a=1$, $r_1 = 12/\sqrt{145}$, and $r_2 = 1/\sqrt{145}$. The model-free GESC converges to a limit cycle as predicted by the average system, whereas the novel ESC visually appears to converge to the origin at this scale. Both model-free ESCs are practically stable to the origin, although the distance from the origin in the limit $t\to\infty$ as a function of $a$ is drastically different.}
    \label{fig:motivating-example:model-free-trajectories}
\end{figure}

Up to this point, we have yet to consider $s$ where $\omega'_2/\omega'_1 \neq 3$. However, one can see the amount of effort required to determine whether the average system is sGPAS to the origin with respect to $a$, since this property needs to be verified for every other possible choice of $(\omega'_1,\omega'_2)$. For $J$ in this example, this task can be done since $J$ is not too complex and the dimensionality is low ($n=2$). For any $\vert\omega'_2/\omega'_1\vert \notin\lbrace 1, 3, \frac{1}{3} \rbrace$, the averaged gradient estimate is the same as \eqref{eq:gesc-example:average-system} except $A$ is now
\begingroup
\renewcommand*{\arraystretch}{2}
\begin{equation*}
    A(r_1, r_2) = \begin{bmatrix}
		6 r_1^2 + 6 r_2^2 & 3 r_1^2 + 6 r_2^2 \\
		6 r_1^2 + 3 r_2^2 & 6 r_1^2 + 3 r_2^2
	\end{bmatrix}
\end{equation*}
\endgroup
with the eigenvalues
\begin{equation*}
    \lambda_{1,2} = 6r_1^2 + \frac{9}{2}r_2^2 \pm \frac{3}{2}\sqrt{8r_1^4 + 20 r_1^2 r_2^2 + 9r_2^4}
\end{equation*}
Both of these eigenvalues are strictly positive for $r_1,r_2\neq 0$, so the average system is stable for almost all choices of the dither amplitude and rate, but not for every choice. Consequently, we are only able to conclude that average systems are, in general, only sGPAS rather than the stronger GAS. Now consider the more expansive task of determining the stability properties over a family of $J$, such as $J$'s which are strictly convex. Rather than going through all these calculations over the entire family of $J$, Theorem \ref{thm:model-based-implies-model-free} provides a sufficient condition to unburden the system analyst from all these considerations. This computational saving is true even when we cannot readily give the average system. While the GESC average system was readily calculated as a polynomial in this example, the average system of \eqref{eq:ex-sys-1} is not easily expressed and can likely only be numerically determined. However, the model-based GESC being GAS is equivalent to \eqref{eq:ex-sys-1:model-based} being GAS, so the model-free ESC in \eqref{eq:ex-sys-1} and its average system are practically stable by Theorem~\ref{thm:model-based-implies-model-free}. This conclusion has been reached without ever calculating the average system of \eqref{eq:ex-sys-1} and is corroborated by the model-free system trajectories shown in Fig.~\ref{fig:motivating-example:model-free-trajectories}.

\section{Proof of Main Results}
    \label{sec:main-results}

The proof of Theorem \ref{thm:model-based-implies-model-free} is broken into two parts, the proof that the average system is sGPAS with respect to the small parameter $a$ and then the proof that the original model-free system is sGPUAS with respect to the small parameter vector $(a, \omega^{-1})$. Both parts will show that each of these systems meets the appropriate definitions of sGPAS and sGPUAS given in Appendix~\ref{sec:practical-stability-definitions}, for which, without loss of generality, we have assumed $\theta^* = 0$ and $x^* = 0$. The motivation behind these proofs is to be as constructive as possible. The system analyst has likely already found a Lyapunov $V$ to prove the model-based ESC is GAS, so it would be efficient to use that same $V$ in order to prove the practical stability of the average and model-free ESC systems. This differs from the proofs for framework-based results, such as \cite{ref:nesic-2010,ref:nesic-2013}, which rely on the existence of a function $\beta\in\mathcal{KL}$ as an outer envelope for the trajectories of the model-based ESC system. Since a system analyst has likely found $V$ but omitted the tedious calculation of $\beta$, we utilize $V$ rather than $\beta$ in the subsequent proofs.

\subsection{The Average System is sGPAS}
	\label{sec:main-results:average-system}

We will use a Lyapunov function to prove that the average system meets the definitions of stability, boundedness, and attractivity. We have assumed that $f$ in \eqref{eq:perturbation-esc} is locally Lipschitz with respect to $\hat{x}$ for the existence of the average system, so $f$ is continuous and $\bar{f}$ is also continuous by \eqref{eq:average-system}. Since $\bar{f}(\cdot, a)$ converges compactly to $h(\cdot)$ as $a\to 0$, $h$ in \eqref{eq:perturbation-free-esc} is continuous. Therefore, trajectories of the model-based ESC exist by Peano's existence theorem although we are not explicitly assuming the uniqueness of these trajectories. Since the model-based system is assumed to be GAS by the theorem statement, and since GAS holds for \emph{every} (potentially nonunique) solution of the model-based system, we can invoke Kurzweil's converse Lyapunov theorem with only the existence of solutions for the model-based ESC. Thus, there exists an autonomous, continuously differentiable $V:\mathbb{R}^m\to\mathbb{R}$, two class $\mathcal{K}_\infty$ functions ($\alpha_1$ and $\alpha_2$) \cite[Ch~4]{ref:khalil-2002}, and a positive definite function $W$ that satisfies the conditions
\begin{gather}
	\alpha_1(\Vert z \Vert) \leq V(z) \leq \alpha_2(\Vert z \Vert) \\
	\label{eq:converse-laypunov-thm:derivative}
	\frac{\partial V}{\partial z}(z) \cdot h(z) \leq -W( z )
\end{gather} 
by the converse Lyapunov theorem \cite[Th~4]{ref:kurzweil-1956}. While there are other converse Lyapunov theorems which may replace $W(z)$ with $\alpha_3(\Vert z \Vert)$ in the right hand side (RHS) of \eqref{eq:converse-laypunov-thm:derivative} \cite[pg~440]{ref:sontag-1989}, the system analyst is more likely to find a $V$ whose Lie derivative satisfies \eqref{eq:converse-laypunov-thm:derivative} when attempting to prove that the model-based ESC is GAS. Therefore, we will proceed with the proof using a Lyapunov function whose Lie derivative satisfies \eqref{eq:converse-laypunov-thm:derivative}.

If we rewrite the vector field $\bar{f}$ as
\begin{equation}
	\bar{f}(\hat{\bar{x}}, a) = h(\hat{\bar{x}}) + R_f(\hat{\bar{x}}, a)
\end{equation}
where $R_f$ is defined as
\begin{equation}
    R_f(\hat{\bar{x}}, a) = \bar{f}(\hat{\bar{x}}, a) - h(\hat{\bar{x}})\,,
\end{equation}
then we can bound the Lie derivative of $V$ along the trajectories of the average system with
\begin{align}
	\frac{\partial V}{\partial \hat{\bar{x}}}(\hat{\bar{x}}) \cdot \rev{\bar{f}}(\hat{\bar{x}}, a) &{}={} \frac{\partial V}{\partial \hat{\bar{x}}} (\hat{\bar{x}}) \cdot h(\hat{\bar{x}}) + \frac{\partial V}{\partial \hat{\bar{x}}} (\hat{\bar{x}}) \cdot R_f(\hat{\bar{x}}, a) 
	\label{eq:lie-derivative:along-average-system:upper-bound}
	{}\leq{} -W(\hat{\bar{x}}) + \Vert R_f(\hat{\bar{x}},a)\Vert \left\Vert \frac{\partial V}{\partial\hat{\bar{x}}}(\hat{\bar{x}}) \right\Vert
\end{align} 
Since $\bar{f}(\hat{\bar{x}}, a) \to h(\hat{\bar{x}})$ compactly as $a\to 0$, $\Vert R_f(\hat{\bar{x}},a) \Vert \to 0$ compactly as $a\to 0$. Thus, for any $\rho_1,\rho_2$ with $0 < \rho_1 < \rho_2 < \infty$ and any margin $c_4\in(0,1)$ fixed once and for all, we can find $a^* = a^*(\rho_1,\rho_2,c_4) > 0$ such that
\begin{equation}
    \label{eq:a-star-margin}
    \Vert R_f(\hat{\bar{x}}, a)\Vert \left\Vert \frac{\partial V}{\partial \hat{\bar{x}}}(\hat{\bar{x}})\right\Vert \leq (1 - c_4)\,W(\hat{\bar{x}})\,,\qquad \forall\, \Vert\hat{\bar{x}}\Vert\in[\rho_1,\rho_2]\,,\ \forall\, a\in(0,a^*)\,,
\end{equation}
which is possible because $\Vert\partial V/\partial\hat{\bar{x}}\Vert$ is bounded and $W$ is bounded away from zero on the compact annulus, while $\Vert R_f(\cdot,a)\Vert$ tends to zero uniformly there. Combining \eqref{eq:a-star-margin} with \eqref{eq:lie-derivative:along-average-system:upper-bound} gives
\begin{equation}
    \label{eq:lie-derivative-margin}
    \frac{\partial V}{\partial \hat{\bar{x}}}(\hat{\bar{x}})\cdot \bar{f}(\hat{\bar{x}},a) \leq -c_4 W(\hat{\bar{x}}) < 0\,,\qquad \forall\, \Vert\hat{\bar{x}}\Vert\in[\rho_1,\rho_2]\,,\ \forall\, a\in(0,a^*)\,.
\end{equation}
The margin $c_4$ renders the decay rate uniform over $a\in(0,a^*)$, which is used in the selection of $c_3$ below. We use $V$ and selection of $\rho_1,\rho_2$ to prove that the average system meets the three required definitions of sGPAS; that the average system is practically stable, practically bounded, and semiglobally practically attractive. 

Note that these proofs are different from those found in \cite{ref:moreau-2000}, which use the assumption of ``convergence of trajectories'' \cite[Hyp~2]{ref:moreau-2000}. The ``convergence of trajectories'' assumption would mean that the trajectories of the average system and model-based systems, starting from the same initial conditions, could be made arbitrarily close over a finite time interval through the choice of ever smaller $a$. A practical stability proof using this assumption would require induction, as was done in \cite{ref:moreau-2000}. The route taken here does not use that hypothesis. Establishing trajectory convergence between the average and model-based systems would call for additional assumptions, such as local Lipschitzness of the model-based system \eqref{eq:perturbation-free-esc} in order to invoke \cite[Th~3.4]{ref:khalil-2002}, and the generality sought in Assumption~\ref{asmp:averaging} (A4) admits model-based vector fields that are merely continuous. Existence of solutions to initial value problems for \eqref{eq:perturbation-free-esc} is all that the argument below requires; uniqueness is not needed.

\subsubsection{The Average System is Practically Stable}

Given some final radius $c_2 > 0$, we select $c_1$, $\rho_1$, and $\rho_2$ such that
$\alpha_2(c_1) < \alpha_1(c_2)$, $0 < \rho_1 < c_1$, and $c_2 < \rho_2$. The boundary of the level set $V_{\alpha_2(c_1)} = \lbrace \hat{\bar{x}}\ \vert\ V(\hat{\bar{x}}) < \alpha_2(c_1) \rbrace$ is within the difference of balls $\mathcal{B}_{\rho_2}\setminus\mathcal{B}_{\rho_1}$ and thus $\frac{d V}{dt} < 0$ everywhere on the boundary whenever $a\in (0,{a^*}(\rho_1,\rho_2))$. Therefore, this level set is forward invariant. From the above selection of the constants $c_1$, $\rho_1$, and $\rho_2$, we can see
\begin{equation}
	\label{eq:average-system-stable-and-bounded}
	\left\Vert\hat{\bar{x}}(t_0)\right\Vert < c_1 \implies \hat{\bar{x}}(t_0)\in V_{\alpha_2 (c_1)} \implies \hat{\bar{x}}(t) \in V_{\alpha_2 (c_1)} \subseteq  V_{\alpha_1 (c_2)} \subseteq \mathcal{B}_{c_2}\ \forall t\in [t_0,\infty)
\end{equation}
when $a \in (0, a^* (\rho_1, \rho_2))$ which is the definition of practical stability that we were trying to achieve. We conclude that the average system is practically stable with respect to the small parameter $a$. All sets and relations can be seen in Fig.~\ref{fig:average-system:stability-and-boundedness}.

\begin{figure}[t!]
	\centering
	\begin{tikzpicture}
	\draw[ultra thick, black, label={$c_1$}] (0,0) circle (0.75cm);
	\begin{scope}[rotate=45]
		\draw[red, ultra thick, label={$V = \alpha_2(c_1)$}] (0,0) ellipse  (0.85 and 1.75);
	\end{scope}
	\draw[ultra thick, black] (0,0) circle (2);
	\draw[dashed] (0,0) circle (2.25);
	\draw[dashed] (0,0) circle (0.5);
	
	\node at (0, .25) {$\rho_1$};
	\node at (0, 2.5) {$\rho_2$};
	\node at (-45:1) {$c_1$};
	\node at (45:1.75) {$c_2$};
	\node[red] at (95: 1.5) {$\partial V_{\alpha_2(c_1)}$};
\end{tikzpicture}
	\caption{Graphical representation of the proof for the average system being practical stability and bounded. Balls with radii $\rho_1$ and $\rho_2$ are shown as dashed while  balls with radii $c_1$ and $c_2$ are solid. $\partial V_{\alpha_2 (c_1)}$ is the boundary of the forward invariant levels set.}
    \label{fig:average-system:stability-and-boundedness}
\end{figure}

\subsubsection{The Average System is Practically Bounded}

The proof of practical boundedness follows in a way similar to the proof of practical stability.  Given some initial radius $c_1 > 0$, we select \rev{$c_2$}, $\rho_1$, and $\rho_2$ to meet once again the requirements of $\alpha_2(c_1) < \alpha_1(c_2)$, $0 < \rho_1 < c_1$, and $c_2 < \rho_2$. With $a\in (0, a^* (\rho_1, \rho_2))$, the level set $V_{\alpha_2 (c_1)}$ is again forward invariant and thus \eqref{eq:average-system-stable-and-bounded} is once again true. We therefore conclude that the average system is practically bounded with respect to the small parameter $a$.

\subsubsection{The Average System is Semiglobally Practically Attractive}

Given an initial and final radius of $c_1 > 0$ and $c_2 > 0$, we select $\rho_1$ so that $0 < \rho_1 < \alpha_{2}^{-1}\circ\alpha_{1}(\min(c_1, c_2))$ and $\rho_2$ so that $\alpha_1^{-1}\circ\alpha_2(\max(c_1, c_2)) < \rho_2$. Note that the previous proofs of practical stability and boundedness explicitly had the ordering $c_1 < c_2$ but the definition of practical attractiveness does not have this ordering. Hence, our selection of $\rho_1$ and $\rho_2$ must be based on the minimum and maximum $c$'s, respectively, but our selection means that the level sets $V_{\alpha_1 (c_2)} = \lbrace \hat{\bar{x}}\ \vert\ V(\hat{\bar{x}}) < \alpha_1 (c_2)\rbrace$ and $V_{\alpha_2 (c_1)} = \lbrace \hat{\bar{x}}\ \vert\ V(\hat{\bar{x}}) < \alpha_2 (c_1)\rbrace$ are forward invariant sets when $a\in (0, a^*(\rho_1, \rho_2))$. Thus, 
\begin{equation}
	\hat{\bar{x}}(t_0)\in V_{\alpha_1 (c_2)} \implies \Vert \hat{\bar{x}}(t) \Vert < c_2\ \forall t \in [t_0,\infty)
\end{equation}
If $\mathcal{B}_{c_1} \subseteq V_{\alpha_1(c_2)}$, then all the starting trajectories are within $V_{\alpha_1 (c_2)}$ and our proof is done with $T(c_1, c_2)\equiv 0$. Otherwise, we must consider the trajectories when $\hat{\bar{x}}(t_0)~\in~ \mathcal{B}_{c_1}~\setminus~V_{\alpha_1(c_2)}$. These trajectories are still in the forward invariant level set $V_{\alpha_2(c_1)}$ by choice of $\rho_1$ and $\rho_2$ so it is now a question of whether there exists a trajectory $\hat{\bar{x}}(t)\in V_{\alpha_2(c_1)}\setminus  V_{\alpha_1(c_2)}$ for all $t\in [t_0,\infty)$. The existence of such a trajectory would serve as a counter-example to the average system being practically attractive.

Note that $V$ would decrease along such a trajectory at least as fast as $c_3$ where 
\begin{equation}
	c_3 = c_4 \min_{\rho_1 \leq \Vert \hat{\bar{x}} \Vert \leq \rho_2} W(\hat{\bar{x}})
\end{equation}
by \eqref{eq:lie-derivative-margin}, with $c_3 > 0$ because $W$ is positive definite and continuous and the annulus is compact. After a time $T(c_1,c_2)$ where 
\begin{equation}
	T(c_1,c_2) = \max\left(\frac{\alpha_2(c_1) - \alpha_1(c_2)}{c_3}, 0\right)
\end{equation}
has elapsed, we can form the upper bound of $V(\hat{\bar{x}}(t_0 + T))$ for trajectories with $\hat{\bar{x}}(t_0) \in V_{\alpha_2(c_1)}\setminus  V_{\alpha_1(c_2)}$ by
\begin{align}
	V(\hat{\bar{x}}(t_0 + T))
        \label{eq:proof:average-system:lyapunov-trajectory-upper-bound}
        &= V(\hat{\bar{x}}(t_0)) + \int_{t_0}^{t_0 + T} \underbrace{\frac{\partial V}{\partial \hat{\bar{x}}}(\hat{\bar{x}}(\tau))^T \bar{f}(\hat{\bar{x}}(\tau), a)}_{\leq -c_3\ \forall \hat{\bar{x}}(t)\in V_{\alpha_2(c_1)}\setminus  V_{\alpha_1(c_2)}} d\tau \\
	& \leq \alpha_2(c_1) - c_3 T \leq \alpha_1(c_2)
\end{align}
This means that $\hat{\bar{x}}(t_0 + T)$ is inside the forward invariant set $V_{\alpha_1(c_2)}$ and thus
\begin{equation}
	\Vert \hat{\bar{x}}(t_0) \Vert < c_1 \implies \hat{\bar{x}}(t) \in V_{\alpha_1(c_2)} \subseteq \mathcal{B}_{c_2}\ \forall t\in[t_0 + T, \infty)
\end{equation}
Since $c_1$ can be arbitrarily large, we conclude that we meet the definition of $\delta$-practically attractive ($\delta$-PA) with $\delta$ being able to be made arbitrarily large by choice of $a$. 

\begin{rem}
	\em
	\label{rem:average-system-convergence-speed}
	The selection of $a^*$ has a tremendous effect on the guaranteed speed of practical convergence of the average system. It is possible to find a sufficiently small $a^*$ for any given $\hat{\bar{x}}$ such that, for all $a\in(0,a^*)$,
	\begin{equation}
		W(\hat{\bar{x}}) - \Vert R_f(\hat{\bar{x}}, a) \Vert \left\Vert \frac{\partial V}{\partial \hat{\bar{x}}}(\hat{\bar{x}}) \right\Vert \geq c_4 W(\hat{\bar{x}})
	\end{equation}
	where $c_4\in (0, 1)$ by $R_f(\hat{\bar{x}}, a) \to 0$ as $a\to 0$, which is \eqref{eq:a-star-margin} written in terms of the decay rate. What this means is that the average system will practically converge to the origin but the convergence speed that we can guarantee is a slower convergence than the model-based system. For the time scaling below to be legitimate, $\beta$ must be the envelope generated from $V$ and $W$ through the Comparison Lemma \cite[Lem~3.4]{ref:khalil-2002} rather than an arbitrary pre-existing $\mathcal{KL}$ bound: replacing $W$ by $c_4 W$ in the comparison differential inequality is what rescales the time argument of its solution by $c_4$, and no such conclusion is available for a $\beta$ obtained by other means. For example, if $\exists\beta\in\mathcal{KL}$ so generated such that
	\begin{equation}
		\Vert z(t_0) \Vert \leq c_1 \implies \Vert z(t) \Vert \leq \beta(\Vert z(t_0) \Vert, t - t_0)
	\end{equation}
	for all $t\in [t_0,\infty)$, then, for every $c_2 > 0$, there exists $a^*(c_1,c_2) = a^*$ such that for all $a\in(0,a^*)$, we have 
	\begin{equation}
		\label{eq:average-system:kl-bound}
		\Vert \hat{\bar{x}}(t_0) \Vert \leq c_1 \implies \Vert \hat{\bar{x}}(t) \Vert \leq \beta(\Vert \hat{\bar{x}}(t_0) \Vert, c_4(t - t_0)) + c_2
	\end{equation}
	for $t\in [t_0,\infty)$. Here, $c_4$ effectively ``slows down'' the guaranteed practical convergence of the average system compared to the model-based system. However, this use of $\beta$ is actually misleading. Our analysis is based on worst-case guarantees, which do not mean that the average system always converges slower than what the model-based system predicts. Consider the GESC acting on the scalar cost function $ J(\theta) = \theta^4 $ where $\theta^* = 0$. The average system is
	\begin{equation}
		\frac{d}{dt} \hat{\bar{\theta}} = -4k\hat{\bar{\theta}}^3 - 3 k a^2 \hat{\bar{\theta}}
	\end{equation}
	Clearly, the additional $- 3 k a^2 \hat{\bar{\theta}}$ aids the average system in converging to $0$, with large $a$ actually leading the average system to converge significantly faster than the model-based system. For this example, we choose $c_4 = 1$ and simply drop the $c_2$ term from \eqref{eq:average-system:kl-bound} to form an outer envelope for the trajectories that is valid for all $a\in(0,\infty)$. One must remember that practical stability proofs are about showing the existence of parameter intervals for which trajectories converge close enough to the equilibrium under a worst-case analysis. Following the steps outlined in our proof, the calculation of $a^*$ for $V(\hat{\bar{\theta}}) = \hat{\bar{\theta}}^2$ yields $a^*(\rho_1,\rho_2) = \sqrt{\frac{4}{3}(1-c_4)}\,\rho_1$; the binding radius is the inner one, since $W = 8k\hat{\bar{\theta}}^4$ while the discarded perturbation term is $6ka^2\hat{\bar{\theta}}^2$, so their ratio is worst at small $\vert\hat{\bar{\theta}}\vert$. Our proof ensures that $\frac{d}{dt}V(\hat{\bar{\theta}}) \leq - 8 \rev{k} c_4 \hat{\bar{\theta}}^4$ for all $\hat{\bar{\theta}}$ and $a$ such that $\rho_1 \leq \vert \hat{\bar{\theta}} \vert \leq \rho_2 $ and $a\in(0,a^*)$, even though the inequality on $\frac{d}{dt}V(\hat{\bar{\theta}})$ holds everywhere; the gap is an artifact of bounding $\Vert R_f\Vert\Vert\partial V/\partial\hat{\bar{\theta}}\Vert$ without regard to the sign of $R_f$, which in this example is stabilizing.
\end{rem}

\begin{rem}
	\em
    Our proof of the average system being practically stable relies on the existence of a differentiable $V$. As mentioned before, a system analyst typically has such a $V$ at hand from proving the model-based system was GAS but perhaps they have only found a $V$ which is continuous and not differentiable. Differentiability of $V$ is an overly restrictive requirement. For example, \cite{ref:wang-2026} considered continuous, but not differentiable, Lyapunov constructions for the model-based system arising from adaptive gradient methods. One such system is
    \begin{subequations}
        \label{eq:example:interconnected-systems:model-based}
    \begin{align}
        \frac{d}{dt} \vartheta &= -\frac{k \nabla J(\vartheta)}{\sqrt{\Vert \varphi \Vert^2 + \eta}} \\
        \frac{d}{dt} \varphi &= \omega_{\ell}\left(\nabla J(\vartheta) - \varphi\right)
    \end{align}
    \end{subequations}
    where $\eta > 0$ is a constant and $J:\mathbb{R^n}\to\mathbb{R}$ is a differentiable cost function that, without loss of generality, is assumed to have a minimum $J(0) = 0$. The Lyapunov function constructed for this system in \cite{ref:wang-2026} was
    \begin{equation}
        \label{eq:example:interconnected-systems:lyapunov-function}
        V(\vartheta, \varphi) = J(\vartheta) + \max\lbrace \Vert \varphi \Vert^2, \gamma\circ J(\vartheta)\rbrace
    \end{equation}
    where $\gamma\in\mathcal{K}$ was determined based on the growth of $\nabla J$ away from the equilibrium. This $V$ is radially unbounded if $J$ is radially unbounded and its upper right Dini derivative of $V$ along the trajectories
    \begin{equation}
        D^+ V \coloneqq \limsup_{h\to 0^+} \frac{1}{h}\left[ V(\vartheta(t + h), \varphi(t + h)) - V(\vartheta(t), \varphi(t)) \right]
    \end{equation}
    was shown to have a negative definite upper bound of
    \begin{equation}
        D^+ V(\vartheta,\varphi) \leq -W(\vartheta,\varphi) \coloneqq \frac{-k \Vert\nabla J(\vartheta)\Vert^2}{\sqrt{\Vert \varphi \Vert^2 + \eta}} + \begin{cases}
            -\omega_{\ell}\Vert \varphi \Vert ^2 & \text{ if } \Vert \varphi \Vert^2 > \gamma\circ J(\vartheta) \\ 
            0 & \text{ otherwise }
        \end{cases}
    \end{equation}
    This lower semicontinuous bound was enough to prove the model-based system was GAS. Pertinent to this work, consider the model-free system
    \begin{subequations}
        \label{eq:example:interconnected-systems:model-free}
    \begin{align}
        \frac{d}{dt}\hat{\theta} &= -\frac{k \hat{g}\left(\omega t, \hat{\theta}, a \right)}{\sqrt{\Vert\hat{\phi}\Vert^2 + \eta}}  \\
        \frac{d}{dt} \hat{\phi} &= \omega_{\ell}\left(\hat{g}\left(\omega t, \hat{\theta}, a \right) - \hat{\phi}\right)
    \end{align}
    \end{subequations}
    whose model-based system is \eqref{eq:example:interconnected-systems:model-based}. Can we prove the average system is practically stable using the $V$ given in \eqref{eq:example:interconnected-systems:lyapunov-function}? The answer is yes: what matters in this case is not the differentiability of $V$ but rather $J$ and $\Vert\varphi\Vert^2$. For any given $0 < \rho_1 < \rho_2 < \infty$, we can partition the annulus 
    $
        \mathcal{D} = \left\lbrace \left(\hat{\bar{\theta}}, \hat{\bar{\phi}}\right)\ \Big\vert\ \rho_1 \leq \left\Vert (\hat{\bar{\theta}}, \hat{\bar{\phi}})\right\Vert \leq \rho_2 \right\rbrace$
    into $\mathcal{D}_{\leq} = \left\lbrace \left(\hat{\bar{\theta}}, \hat{\bar{\phi}}\right)\in\mathcal{D}\ \vert\ \Vert\hat{\bar{\phi}}\Vert^2 \leq \gamma\circ J(\hat{\bar{\theta}}) \right\rbrace$ and $\mathcal{D}_{\geq} = \left\lbrace \left(\hat{\bar{\theta}}, \hat{\bar{\phi}}\right)\in\mathcal{D}\ \vert\ \Vert\hat{\bar{\phi}}\Vert^2 \geq \gamma\circ J(\hat{\bar{\theta}}) \right\rbrace$.
    Then, for any $c_4\in (0,1)$, we can always find an $a^*$ such that, for all $a\in(0,a^*)$,
    \begin{equation*}
            \left\Vert \frac{\nabla J(\hat{\bar{\theta}}) \cdot R_g(\hat{\bar{\theta}}, a)}{\sqrt{\Vert\hat{\bar{\phi}}\Vert^2 + \eta}} \right\Vert \leq (1-c_4) \frac{\Vert \nabla J(\hat{\bar{\theta}})\Vert^2}{\sqrt{\Vert\hat{\bar{\phi}}\Vert^2 + \eta}} \qquad \forall\ (\hat{\bar{\theta}}, \hat{\bar{\phi}}) \in \mathcal{D}_{\leq}
    \end{equation*}
    and
    \begin{equation*}
            \left\Vert \frac{\nabla J(\hat{\bar{\theta}}) \cdot R_g(\hat{\bar{\theta}}, a)}{\sqrt{\Vert\hat{\bar{\phi}}\Vert^2 + \eta}} + 2\omega_{\ell} \hat{\bar{\phi}}\cdot R_g(\hat{\bar{\theta}}, a) \right\Vert \leq (1-c_4) \left[\frac{\Vert \nabla J(\hat{\bar{\theta}})\Vert^2}{\sqrt{\Vert\hat{\bar{\phi}}\Vert^2 + \eta}} + \omega_{\ell}\Vert\hat{\bar{\phi}}\Vert^2 \right] \qquad \forall\ (\hat{\bar{\theta}}, \hat{\bar{\phi}}) \in \mathcal{D}_{\geq}
    \end{equation*}
    by $R_g(\cdot, a) \to 0$ compactly in the limit $a\to 0$. Thus, the upper right Dini derivative of $V$ along trajectories in the average system is bounded by
    \begin{equation}
        D^+ V\left(\hat{\bar{\theta}}, \hat{\bar{\phi}}\right) \leq - c_4 W \left(\hat{\bar{\theta}}, \hat{\bar{\phi}}\right) \quad \forall \left(\hat{\bar{\theta}}, \hat{\bar{\phi}}\right) \in \mathcal{D} \text{ and } a\in(0,a^*)
    \end{equation}
    This negative definite bound in the annulus $\mathcal{D}$ is enough to establish that any level set of $V$ contained within $\mathcal{D}$ is forward invariant for the average system with any $a\in(0,a^*)$. Additionally, we can use the Comparison Lemma \cite[Lem~3.4]{ref:khalil-2002} to guarantee $V$ is decreasing for all average system trajectories by at least
    \begin{equation*}
        c_3 = \min_{\rho_1 \leq \Vert (\hat{\bar{\theta}}, \hat{\bar{\phi}}) \Vert\leq \rho_2} c_4 W(\hat{\bar{\theta}}, \hat{\bar{\phi}})
    \end{equation*}
    as long as the trajectories are within the annulus. The forward invariance and the guaranteed minimum decay rate of $V$ are enough to substitute this continuous $V$ into the steps given in Section~\ref{sec:main-results:average-system}. The only difference is that we use a differential inequality in \eqref{eq:proof:average-system:lyapunov-trajectory-upper-bound} and the comparison lemma to establish that the average system is semiglobally practically attractive in the proof by contradiction.
\end{rem}

\subsection{The Model-Free ESC is sGPUAS}
	\label{sec:main-results:model-free}

The analysis of the model-free ESC system will show that the system is sGPUAS to the origin by proving the system meets the required definitions. This analysis is aided by the fact that the average system is sGPAS and that we can make trajectories of the model-free and average systems arbitrarily close from the same initial condition, for any given $a$, over a finite time interval by taking $\omega$ arbitrarily large \cite[Th~4.3.6]{ref:sanders-2007}, \cite[Th~10.5]{ref:khalil-2002}. The theorems for this arbitrary closeness are generally given using the ``fast'' time scale $\tau=\omega t$ with the derivative $\frac{d}{d\tau} = \omega^{-1}\frac{d}{dt} = \varepsilon\frac{d}{dt}$ which generates the differential equations
\begin{align}
	\frac{d}{d\tau} \hat{x} &= \varepsilon f(\tau, \hat{x}, a) \\
	\frac{d}{d\tau} \hat{\bar{x}} &= \varepsilon \bar{f}(\hat{\bar{x}}, a)
\end{align}
The closeness for general averaging is then given in terms of $\mathcal{O}(\alpha(\varepsilon))$ for some $\alpha\in\mathcal{K}$ over the time interval $\tau\in[\tau_0, \tau_0 + T/\varepsilon]$. This is why we can select a larger $\omega$ (equivalent to a smaller $\varepsilon$) to force the trajectories between the systems to be arbitrarily close over the ``slow'' time interval $t\in [t_0, t_0 + T]$.

However, because these are singular perturbation theorems, one must consider how values of $a$ would affect the choice of $\omega$. As $a$ likely influences the local Lipschitz constant of $f$ with respect to $\hat{x}$ in \eqref{eq:perturbation-esc}, $a$ impacts the growth of the bounds between trajectories of the two systems, as seen in the details of the proof of \cite[Lem~4.2.7]{ref:sanders-2007}. Thus, we summarize the ability to make the trajectories arbitrarily close in the context of this work with the following lemma.

\begin{lem}[Adapted from {\cite[Th~10.5]{ref:khalil-2002}}]
    \label{lem:averaging-closeness}
	Consider trajectories for systems \eqref{eq:perturbation-esc} and \eqref{eq:average-system} and compact sets $\mathcal{D}_0 \subset \mathcal{D}\subset\mathbb{R}^m$ such that for the initial conditions $\hat{\bar{x}}(t_0)\in\mathcal{D}_0$, the trajectories satisfy $\hat{\bar{x}}(t)\in\mathcal{D}$ for all $t\in[t_0,t_0 + T]$ and some $T \in (0,\infty)$. Let Assumption~\ref{asmp:averaging} hold. Then, for a given $a > 0$ and $\Delta > 0$, $\exists\omega^* = \omega^*(T, \Delta, a) > 0$, independent of $t_0$, such that if $\hat{x}(t_0) = \hat{\bar{x}}(t_0)\in\mathcal{D}_0$, then $\forall \omega \in( \omega^*,\infty)$, the trajectories of the original and average system satisfy
		\begin{equation}
			\left\Vert \hat{x}(t) - \hat{\bar{x}}(t) \right\Vert < \Delta,\ \forall t\in [t_0, t_0 + T]
	\end{equation}
\end{lem}

The independence of $\omega^*$ from $t_0$ is used in what follows, since the induction used in the subsequent proofs restarts the comparison on each interval $[t_n, t_{n+1}]$, at which point the fast phase of the nonautonomous system \eqref{eq:perturbation-esc} is $\omega t_n$ rather than $\omega t_0$. Autonomy of the average system says nothing about that phase. What supplies the required uniformity is (A2): the error committed by replacing $f$ with $\bar{f}$ over a window is bounded by $\kappa_{\mathcal{D}}(a)\sigma(T)$ \emph{for every} starting phase $\tau_0$, and the Gr\"onwall estimate underlying \cite[Th~10.5]{ref:khalil-2002} then produces a bound depending only on $T$, $\Delta$, $a$, and the Lipschitz and boundedness constants of (A1) on $\mathcal{D}$. For a periodic dither this is immediate from compactness of the phase; for an almost periodic dither it is the content of \eqref{eq:strong-average}.

The remainder of this section merely verifies that the original model-free ESC system meets the required definitions of sGPUAS; that the model-free ESC system is practically stable, practically bounded, and semiglobally practically uniformly attractive. The proofs will follow similarly to the proof of \cite[Th~1]{ref:moreau-2000} but with average system being sGPAS rather than GAS. The hypotheses here differ from those of \cite[Lem~1]{ref:tan-2005}, which rests on \cite[Th~1]{ref:nesic-2001}: the latter is stated for an average system that is globally practically asymptotically stable (GPAS), whereas Section~\ref{sec:main-results:average-system} establishes the semiglobal property, and it carries a local Lipschitz requirement in $z$ that Assumption~\ref{asmp:averaging} does not impose on the model-based ESC. The proof below therefore proceeds directly from the sGPAS property and Lemma~\ref{lem:averaging-closeness}.

\subsubsection{The Model-Free ESC is Practically Stable}

Given a final radius of $c_2 > 0$, we select parameters $c_1$, $d$, $\rho_1$, and $\rho_2$ such that $\alpha_2(c_1)<\alpha_1(c_2)$, $\alpha_2(d)<\alpha_1(c_1)$, $0<\rho_1<\alpha_{2}^{-1}(d)$, and $\alpha_1(\rho_2) > \alpha_2(c_2)$. By the average system being sGPAS, $\exists\ T(c_1,d) = T,\ {a^*} > 0$ such that $\forall a\in (0,a^*)$
\begin{align}
	\label{eq:model-free:ps:average-system:plaa}
	&\Vert\hat{\bar{x}}(t_0)\Vert < c_1\implies \Vert\hat{\bar{x}}(t_0 + T)\Vert < d \\
	\label{eq:model-free:ps:average-system:pb}
	&\Vert\hat{\bar{{x}}}(t_0)\Vert < c_1\implies \Vert\hat{\bar{{x}}}(t)\Vert < c_2,\ \forall\ t\in [t_0,t_0+T]
\end{align}
 where \eqref{eq:model-free:ps:average-system:pb} is a consequence of the level set $V_{\alpha_2(c_1)} = \lbrace \hat{\bar{x}}\ \vert\ V(\hat{\bar{x}}) < \alpha_2(c_1) \rbrace\subseteq\mathcal{B}_{c_2}$ being forward invariant by choices of $\rho_1$ and $\rho_2$. Choose $\Delta = \min\lbrace c_1 - d, c_2 - \alpha_1^{-1}\circ\alpha_2 (c_1) \rbrace$ so that $\omega^* = {\omega^{*}}(T, \Delta, a)$ and $\forall\ \omega\in ({\omega^*}, \infty)$ 
\begin{equation}
	\label{eq:model-free:ps:closeness}
	\Vert \hat{x}(t) \Vert \leq \Vert \hat{\bar{x}}(t) \Vert + \Delta\ \forall\ t\in [t_0, t_0 + T]
\end{equation}
if $\hat{x}(t_0) = \hat{\bar{x}}(t_0)$ by the closeness of trajectories given in Lemma~\ref{lem:averaging-closeness}. Combining \eqref{eq:model-free:ps:average-system:plaa} and \eqref{eq:model-free:ps:average-system:pb} with \eqref{eq:model-free:ps:closeness} results in the statements
\begin{align}
	\label{eq:model-free:ps:model-free:plaa}
	&\Vert\hat{x}(t_0)\Vert < c_1\implies \Vert\hat{x}(t_0 + T)\Vert < c_1 \\
	\label{eq:model-free:ps:model-free:pb}
	&\Vert\hat{x}(t_0)\Vert < c_1\implies \Vert\hat{x}(t)\Vert < c_2,\ \forall\ t\in [t_0,t_0+T]
\end{align}
where \eqref{eq:model-free:ps:model-free:plaa} is a result of $\Delta \leq c_1 - d$ and \eqref{eq:model-free:ps:model-free:pb} is a result of $\Delta \leq c_2 - \alpha_{1}^{-1}\circ \alpha_2 (c_1)$.

Since the average system is autonomous, we can use the conditions \eqref{eq:model-free:ps:model-free:plaa} and \eqref{eq:model-free:ps:model-free:pb} in an inductive manner by considering the following time intervals $[t_n,t_{n+1}]$ where $t_n = t_0 + nT$, $\Vert \hat{x}(t_n) \Vert < c_1$, and $n\in\lbrace 1,2,\ldots\rbrace$. In these time intervals, we have the statements
\begin{align}
	\label{eq:model-free:ps:model-free:plaa:inductive}
	&\Vert\hat{x}(t_n)\Vert < c_1\implies \Vert\hat{x}(t_{n+1})\Vert < c_1 \\
	\label{eq:model-free:ps:model-free:pb:inductive}
	&\Vert\hat{x}(t_n)\Vert < c_1\implies \Vert\hat{x}(t)\Vert < c_2,\ \forall\ t\in [t_n,t_{n+1}]
\end{align}
due to the average system being autonomous. Thus, we can extend the upper bound of the trajectories over finite time intervals to
\begin{equation}
	\Vert\hat{x}(t_0)\Vert < c_1\implies \Vert\hat{x}(t)\Vert < c_2,\ \forall\ t\in [t_0, \infty)
\end{equation}
and conclude that the model-free ESC is practically stable with respect to the small parameter vector $(a,\omega^{-1})$. See Fig.~\ref{fig:model-free:stability-and-boundedness} for a graphical representation of this proof.

\begin{figure}[t!]
	\centering
	\begin{tikzpicture}
	\draw[dashed] (4, 0) node[anchor=north] {$t_{4}$} -- ++ (0, 3);
	\draw[dashed] (0,1.5) --++(4.75, 0) node[anchor=west] {$c_1$};
	\draw[dashed] (0,0.5) --++(4.75, 0) node[anchor=west] {$d$};
	\draw[dashed] (0,0.5 + 0.65) --++(4.75, 0);
	\draw[latex-latex] (4.25, 0.5) -- node[anchor=west] {$\Delta$} ++(0, 0.65); 
	\draw[dashed] (0,2.6) --++(4.75, 0) node[anchor=west] {$c_2$};
	\draw[dashed] (0,2.6 - 0.65) --++(4.75, 0) node[anchor=west] {$c_2 - \alpha_1^{-1}\alpha_2(c_1)$};
	\draw[latex-latex] (4.25, 2.6) -- node[anchor=west] {$\Delta$} ++(0, -0.65);
	\draw[ultra thick, black, -latex] (0,0) --++(4.5, 0) node[anchor=west] {$t$};
	\draw[ultra thick, black, -latex] (0,0) --++(0, 3.25) node[anchor=south] {\textcolor{blue}{$\Vert \hat{\bar{x}} \Vert$} and \textcolor{red}{$\Vert \hat{x}\Vert$}};
	\foreach \n in {0,...,3}{
		\draw[dashed] (\n, 0) node[anchor=north] {$t_{\n}$} -- ++ (0, 3);
		\draw[ultra thick, blue, {Circle[width=5pt, length=5pt]}-{Circle[open, width=5pt, length=5pt]}, shorten <=-0.1cm, shorten >=-0.1cm]  ({\n}, 1.5) .. controls (\n + 0.5, 1.75) .. (\n + 1, 0.5);
		\begin{scope}[yshift=0.65cm]
			\draw[ultra thick, dashed, blue, {Circle[width=5pt, length=5pt]}-{Circle[open, width=5pt, length=5pt]}, shorten <=-0.1cm, shorten >=-0.1cm]  ({\n}, 1.5) .. controls (\n + 0.5, 1.75) .. (\n + 1, 0.5);
		\end{scope}
		\draw[ultra thick, red, {Circle[width=5pt, length=5pt]}-{Circle[open, width=5pt, length=5pt]}, shorten <=-0.1cm, shorten >=-0.1cm] (\n,1.5) .. controls (\n + 0.5, 2.) .. (\n + 1, 0.9);
	}
\end{tikzpicture}
	\caption{Graphical representation of the proof of the model-free ESC being practically stable and bounded. Bounds for \textcolor{blue}{$\Vert \hat{\bar{x}}\Vert$} are shown with a blue line {\raisebox{2pt}{\protect\tikz \protect\draw[blue, ultra thick] (0,0) --++ (0.5, 0);}}, bounds for \textcolor{red}{$\Vert \hat{x}\Vert$} are shown with a red line {\raisebox{2pt}{\protect\tikz \protect\draw[red, ultra thick] (0,0) --++ (0.5, 0);}}, and $\Vert\hat{\bar{x}}\Vert~+~\Delta$ is shown with a dashed blue line {\raisebox{2pt}{\protect\tikz \protect\draw[blue, ultra thick, dashed] (0,0) --++ (0.5, 0);}}. All curves are inclusive on the left endpoint and exclusive on the right endpoint. These discontinuities are a result of the use of induction in the proof. The dashed blue line forms an upper bound for $\Vert \hat{x}(t) \Vert$ to show that the model-free ESC system meets the definitions of practical stability and boundedness.}
    \label{fig:model-free:stability-and-boundedness}
\end{figure}

\subsubsection{The Model-Free ESC is Practically Bounded}

Given some initial radius $c_1$, we choose the strictly positive constants $c_2$, $d$, $\rho_1$, and $\rho_2$ such that $c_2 > \alpha_1^{-1}\circ\alpha_2(c_1)$, $d<\alpha_2^{-1}\circ\alpha_1(c_1)$, $\rho_1<\alpha_{2}^{-1}(d)$, and $\rho_2 > \alpha_1^{-1}\circ\alpha_2(c_2)$. Then the arguments are identical to the previous practical stability arguments. Thus, the model-free ESC system is practically bounded with respect to the small parameter vector $(a,\omega^{-1})$.

\subsubsection{The Model-Free ESC is Semiglobally Practically Uniformly Attractive}

The proof of the model-free ESC being $\delta$-PUA with $\delta$ being arbitrarily large is broken up into two parts. The initial part is to show that for the initial conditions $\Vert \hat{x}(t_0) \Vert < c_1$, there are critical parameter values $a^*$ and $\omega^*(a)$ such that all those trajectories are within a ball of radius $d_1$ after time $T$ has elapsed when $a\in (0, a^*)$ and $\omega\in (\omega^*(a),\infty)$. The inductive part is to show that since all trajectories are in $\mathcal{B}_{d_1}$ after $T$ seconds have elapsed, they will never leave $\mathcal{B}_{c_2}$.

\paragraph{Initial Part}

Given an initial and final radius of $c_1$ and $c_2$, we choose $d_1,d_2> 0$ such that $d_1 < \alpha_2^{-1}\circ\alpha_1(\min(c_1, c_2))$ and $d_2 < \alpha_2^{-1}\circ\alpha_1(d_1)$. In addition, we must choose $\rho_1$ and $\rho_2$ such that $0 < \rho_1 < \alpha_2^{-1}\circ\alpha_1(d_2)$ and $\rho_2 > \alpha_1^{-1}\circ\alpha_2(\max(c_1,c_2))$.  From the average system being sGPAS, $\exists\ T(c_1, d_2)=T,{a^*}(\rho_1,\rho_2) = a^*$ such that $\forall\ a\in (0,{a^*})$
\begin{equation}
	\Vert\hat{\bar{x}}(t_0)\Vert < c_1 \implies \Vert \hat{\bar{x}}(t_0 + T) \Vert < d_2 < c_2
\end{equation}
Next we make the trajectories between the average and model-free systems to be $\Delta = \min\lbrace c_2 - \alpha_1^{-1}\circ\alpha_2(d_1),\ d_1 - d_2\rbrace $ close over a time interval of $T$ length from the same initial condition using Lemma~\ref{lem:averaging-closeness}. With ${\omega^*} = {\omega^*}(T, \Delta, a)$, $\forall\ \omega\in ({\omega^*}, \infty)$
\begin{equation}
	\Vert\hat{x}(t_0)\Vert < c_1 \implies \Vert \hat{x}(t_0 + T) \Vert < d_1 < c_2
\end{equation}

\paragraph{Inductive Part}

For this inductive part, the analysis will take place over time intervals of $[t_n,t_{n+1}]$ where $t_n = t_0 + nT$ with $n \in \lbrace 1,2,\ldots\rbrace$. From the initial part, and by the average system being autonomous, we know that the same ${a^*}$ allows for
\begin{equation}
	\Vert \hat{\bar{x}}(t_n)\Vert < d_1 < c_1 \implies \Vert \hat{\bar{x}}(t_{n+1})\Vert < d_2
\end{equation}
and
\begin{equation}
	\Vert \hat{\bar{x}}(t_n)\Vert < d_1 \implies \Vert \hat{\bar{x}}(t)\Vert < \alpha_1^{-1}\circ\alpha_2(d_1),\ \forall\ t\in [t_n, t_{n+1}]
\end{equation}
whenever $a\in(0, a^*)$. With the ${\omega^*}$ from the initial part granting the $\Delta$ closeness between trajectories, we can conclude that
\begin{align}
	& \Vert \hat{x}(t_n)\Vert < d_1 \implies \Vert \hat{x}(t_{n+1})\Vert < d_1\\
	& \Vert \hat{x}(t_n)\Vert < d_1 \implies \Vert \hat{x}(t)\Vert < c_2,\ \forall\ t\in [t_n, t_{n+1}]
\end{align}
$\forall\ \omega\in ({\omega^*},\infty)$. With this we can say that $\forall\ a\in (0,{a^*})$ and $\forall\ \omega\in ({\omega^*},\infty)$, the trajectories satisfy $\Vert \hat{x}(t)\Vert < c_2$ $\forall\ t\in [t_0 + T,\infty)$, and thus complete the induction part of the proof. Now we can combine both the initial and inductive parts to state
\begin{equation}
    \Vert \hat{x}(t_0) \Vert < c_1 \implies \Vert \hat{x}(t) \Vert < c_2 \forall\ t\in[t_0+ T,\infty)
\end{equation}
for all $a\in(0,a^*)$ and for all $\omega\in(\omega^*(T, \Delta, a),\infty)$. Hence, the model-free ESC system is $\delta$-PUA with respect to the small parameter vector $(a, \omega^{-1})$ and $\delta$ being able to be made arbitrarily large. A graphical representation of this proof can be found in Fig.~\ref{fig:model-free:d-puaa}.

\begin{figure}[t!]
	\centering
	\begin{tikzpicture}
	\draw[dashed] (6, 0) node[anchor=north] {$t_{6}$} -- ++ (0, 3);
	\draw[dashed] (0,1) --++(6.75, 0) node[anchor=west] {$d_1$};
	\draw[dashed] (0,0.5) --++(6.75, 0) node[anchor=west] {$d_2$};
	\draw[latex-latex] (6.25, 0.5) -- node[anchor=west] {$\Delta$} ++ (0, 0.5);
	\draw[dashed] (0,1.65) --++(6.75, 0) node[anchor=west] {$c_2$};
	\draw[dashed] (0,1.65 - 0.5) --++(6.5, 0);
	\draw[latex-latex] (6.25, 1.65) -- node[anchor=west] {$\Delta$} ++ (0, -0.5);
	\draw[dashed] (0,2) --++(6.75, 0) node[anchor=west] {$c_1$};
	\draw[ultra thick, black, -latex] (0,0) --++(6.5, 0) node[anchor=west] {$t$};
	\draw[ultra thick, black, -latex] (0,0) --++(0, 3) node[anchor=south] {\textcolor{blue}{$\Vert \hat{\bar{x}} \Vert$} and \textcolor{red}{$\Vert \hat{x}\Vert$}};
	\draw[ultra thick, blue, {Circle[width=5pt, length=5pt]}-{Circle[open, width=5pt, length=5pt]}, shorten <=-0.1cm, shorten >=-0.1cm]  (0, 2) .. controls (0.5, 1.25) .. (1, 0.5);
	\begin{scope}[yshift=0.5cm]
		\draw[ultra thick, blue, dashed, {Circle[width=5pt, length=5pt]}-{Circle[open, width=5pt, length=5pt]}, shorten <=-0.1cm, shorten >=-0.1cm]  (0, 2) .. controls (0.5, 1.25) .. (1, 0.5);
	\end{scope}
	\draw[ultra thick, red, {Circle[width=5pt, length=5pt]}-{Circle[open, width=5pt, length=5pt]}, shorten <=-0.1cm, shorten >=-0.1cm] (0,2) .. controls (0.5, 1.35) .. (1, 0.75);
	\draw[dashed] (0, 0) node[anchor=north] {$t_{0}$} -- ++ (0, 3);
	\foreach \n in {1,...,5}{
		\draw[dashed] (\n, 0) node[anchor=north] {$t_{\n}$} -- ++ (0, 2.5);
		\draw[ultra thick, blue, {Circle[width=5pt, length=5pt]}-{Circle[open, width=5pt, length=5pt]}, shorten <=-0.1cm, shorten >=-0.1cm]  ({\n}, 1) .. controls (\n + 0.5, 1.15) .. (\n + 1, 0.25);
		\begin{scope}[yshift=0.5cm]
			\draw[ultra thick, blue, dashed, {Circle[width=5pt, length=5pt]}-{Circle[open, width=5pt, length=5pt]}, shorten <=-0.1cm, shorten >=-0.1cm]  ({\n}, 1) .. controls (\n + 0.5, 1.15) .. (\n + 1, 0.25);
		\end{scope}
		\draw[ultra thick, red, {Circle[width=5pt, length=5pt]}-{Circle[open, width=5pt, length=5pt]}, shorten <=-0.1cm, shorten >=-0.1cm] (\n,1) .. controls (\n + 0.5, 1.5) .. (\n + 1, 0.55);
	}
\end{tikzpicture}
	\caption{Graphical representation of the practical uniform attractivity of the model-free ESC systems.  Bounds for \textcolor{blue}{$\Vert \hat{\bar{x}}\Vert$} are shown with a blue line {\raisebox{2pt}{\protect\tikz \protect\draw[blue, ultra thick] (0,0) --++ (0.5, 0);}}, bounds for \textcolor{red}{$\Vert \hat{x}\Vert$} are shown with a red line {\raisebox{2pt}{\protect\tikz \protect\draw[red, ultra thick] (0,0) --++ (0.5, 0);}}, and $\Vert\hat{\bar{x}}\Vert~+~\Delta$ is shown with a dashed blue line {\raisebox{2pt}{\protect\tikz \protect\draw[blue, ultra thick, dashed] (0,0) --++ (0.5, 0);}}. Again, all curves are inclusive on the left endpoint and exclusive on the right endpoint with discontinuities a result of the use of induction in the proof. The dashed blue line forms an upper bound for $\Vert \hat{x}(t) \Vert$ to show that the model-free ESC system meets the definitions of practical uniform attractiveness.}
    \label{fig:model-free:d-puaa}
\end{figure}

\section{Conclusions}
This work proves that the model-based ESC being GAS is sufficient to prove that the average system is sGPAS and the model-free ESC is sGPUAS. The result follows from relatively weak assumptions which are generally satisfied for most model-free ESC constructions. Importantly, we have relaxed the required assumptions on the model-based ESC; the model-based ESC only needs to have a continuous vector field, not necessarily a locally Lipschitz vector field. This relaxation precludes the approaches used in the previous literature, since a ``closeness of trajectories'' cannot be established between the average and model-based ESC systems due to potential non-uniqueness of trajectories for the model-based ESC system. Instead, we used the compact convergence of the average system to the model-based ESC to form a Lyapunov-based proof. Since the system analyst can design a model-based ESC to ensure the existence of the average system and its compact convergence to the model-based ESC system, they can conclude practical stability results for the model-free ESC by studying the stability of the model-based ESC system without needing to establish the stability results of each average system.

\appendix
\section{Practical Stability Definitions}
	\label{sec:practical-stability-definitions}

Consider the stability of dynamical systems described by the differential equation
\begin{equation}
    \label{eq:definitions:parameter-dependent-ode}
	\frac{dx}{dt} = f(t, x, \varepsilon_1,\ldots,\varepsilon_m)
\end{equation}
The stability of such systems is often given in terms of practical stability \cite{ref:moreau-2000,ref:teel-1998,ref:tan-2006}. The practical stability definitions in the $\epsilon-\delta$ forms are used in the proofs of Section \ref{sec:main-results} and are given below.

\begin{defn}[Practical Stability]	
	\label{def:ps}
	The origin of the system is said to be practically stable (PS) with respect to the small parameter vector $(\varepsilon_1,\ldots,\varepsilon_m)$ if $\forall\ c_2 > 0$, $\exists c_1 > 0$ and $\exists \varepsilon_1^*(c_1, c_2) > 0$ such that $\forall\varepsilon_1\in(0, \varepsilon_1^*)$, $\exists\varepsilon_2^*(\varepsilon_1, c_1, c_2) > 0$ such that $\forall\varepsilon_2\in(0, \varepsilon_2^*),\ldots,\exists \varepsilon_m^*(\varepsilon_1,\ldots,\varepsilon_{m-1}, c_1, c_2)$ such that $\forall \varepsilon_m\in(0,\varepsilon_m^*)$,
	\begin{equation}
		\Vert x(t_0) \Vert < c_1 \implies \Vert x(t) \Vert < c_2 \ \forall\ t\in[t_0,\infty)
	\end{equation}
\end{defn}

\begin{defn}[Practical Boundedness]
	\label{def:pb}
	The origin of the system is said to be practically bounded (PB) with respect to the small parameter vector $(\varepsilon_1,\ldots,\varepsilon_m)$ if $\forall\ c_1 > 0$, $\exists c_2 > 0$ and $\exists \varepsilon_1^*(c_1, c_2) > 0$ such that $\forall\varepsilon_1\in(0, \varepsilon_1^*)$, $\exists\varepsilon_2^*(\varepsilon_1, c_1, c_2) > 0$ such that $\forall\varepsilon_2\in(0, \varepsilon_2^*),\ldots,\exists \varepsilon_m^*(\varepsilon_1,\ldots,\varepsilon_{m-1}, c_1, c_2)$ such that $\forall \varepsilon_m\in(0,\varepsilon_m^*)$,
	\begin{equation}
		\Vert x(t_0) \Vert < c_1 \implies \Vert x(t) \Vert < c_2 \ \forall\ t\in[t_0,\infty)
	\end{equation}
\end{defn}

\begin{defn}[$\delta$-Practical Uniform Attractivity]
	\label{def:dpuaa}
	The origin of the system is said to be $\delta$-practically uniformly attractive ($\delta$-PUA) with respect to the parameter vector $(\varepsilon_1,\ldots,\varepsilon_m)$ if $\forall\ c_1 \in (0,\delta)$ and $c_2 > 0$, $\exists T(c_1, c_2) = T \geq 0$ and $\exists \varepsilon_1^*(c_1, c_2) > 0$ such that $\forall\varepsilon_1\in(0, \varepsilon_1^*)$, $\exists\varepsilon_2^*(\varepsilon_1, c_1, c_2) > 0$ such that $\forall\varepsilon_2\in(0, \varepsilon_2^*),\ldots,\exists \varepsilon_m^*(\varepsilon_1,\ldots,\varepsilon_{m-1}, c_1, c_2)$ such that $\forall \varepsilon_m\in(0,\varepsilon_m^*)$,
	\begin{equation}
		\Vert x(t_0) \Vert < c_1 \implies \Vert x(t) \Vert < c_2 \ \forall\ t\in[t_0 + T,\infty)
	\end{equation}
	If $\delta$ can be made arbitrarily large, then the system is semiglobally practically uniformly attractive (sGPUA). If the system \eqref{eq:definitions:parameter-dependent-ode} is autonomous, then we neglect to include the qualifier ``uniformly" to the attractiveness property since it is redundant for finite-dimensional autonomous systems.
\end{defn}

\begin{defn}[sGPUAS]
	\label{def:sgpuas}
	The origin of the system is said to be semiglobally practically uniformly asymptotically stable (sGPUAS) with respect to the small parameter vector $(\varepsilon_1,\ldots,\varepsilon_m)$ if it is PS, PB, and sGPUA by the small parameter vector $(\varepsilon_1,\ldots,\varepsilon_m)$.
\end{defn}

Alternatively, we combine the aforementioned definitions into an equivalent definition based on a $\mathcal{KL}$ function. This alternative definition is sometimes more convenient to use.

\begin{defn}[{\cite[Def~1]{ref:tan-2006}}]
	\label{def:sgpuas:alternative-definition}
	The system is said to be sGPUAS to the origin with respect to the parameter vector $(\varepsilon_1,\ldots,\varepsilon_m)$ if $\exists\beta\in\mathcal{KL}$ such that for any choice of $c_1 > 0$ and $c_2> 0$, there $\exists \varepsilon_1^*(c_1,c_2) > 0$ such that $\forall\varepsilon_1\in(0, \varepsilon_1^*)$, $\exists\varepsilon_2^*(\varepsilon_1, c_1, c_2)>0$ such that $\forall\varepsilon_2\in(0, \varepsilon_2^*),\ldots,\exists\varepsilon_m^*(\varepsilon_1,\ldots,\varepsilon_{m-1},c_1,c_2)$ such that $\forall\varepsilon_m\in(0,\varepsilon_m^*)$,
	\begin{equation}
		\left\Vert x(t_0) \right\Vert \leq c_1 \implies \left\Vert x(t) \right\Vert \leq \beta\left(\left\Vert x(t_0) \right\Vert, t - t_0\right) + c_2
	\end{equation}
	for all $t\in[t_0,\infty)$.
\end{defn}

If we can remove the dependence on $c_1$ when determining every threshold $\varepsilon_1^*,\ldots,\varepsilon_m^*$, then these definitions are global rather than semiglobal. If we can remove the dependence on $c_2$ when determining the aforementioned thresholds, then we drop the notion of practical stability, and stability follows the classical stability definitions found in \cite[Ch~4]{ref:khalil-2002}. It should be noted that we take the convention that if solutions are non-unique for the initial value problem, then the stability definitions hold over all possible solutions.

\section{Comparison Function Based Proof of Theorem~\ref{thm:model-based-implies-model-free}}
    \label{sec:kl-based-proofs}
    
For readers who prefer the comparison-based definition of practical stability given in Definition~\ref{def:sgpuas:alternative-definition}, and for general completeness, we present alternative proofs using comparison functions.

\subsection{The Average System is sGPAS}
\begin{pf}
Given a $c_1 > 0$ and $c_2 > 0$ for performance and $\beta$ for the original system, define
\begin{equation}
	\label{eq:alpha-beta}
	\alpha_{\beta}(r) := \beta(r, 0)
\end{equation}
Determine a $c_4\in(0,1)$ such that
\begin{equation}
	\label{eq:beta-closeness}
	\sup_{t\in[t_0,\infty)} \max_{r\in[0,c_1]} \left[\beta(r, c_4(t - t_0)) - \beta(r, t - t_0)\right] \leq \frac{c_2}{2},
\end{equation} 
where the existence of $c_4$ is guaranteed by $\beta\in\mathcal{KL}$, specifically by $\beta(r,\cdot)\in\mathcal{L}$. The argument is short enough to give. Fix $S>0$ with $\beta(c_1,S)\leq c_2/2$ and restrict attention to $c_4\in[1/2,1)$. For $t - t_0 \geq S/c_4$ the first term already satisfies $\beta(r, c_4(t-t_0))\leq \beta(c_1,S)\leq c_2/2$ and the bracket in \eqref{eq:beta-closeness} is at most $c_2/2$, since the subtracted term is nonnegative. For $t - t_0 \leq S/c_4 \leq 2S$, the function $\beta$ is uniformly continuous on the compact set $[0,c_1]\times[0,2S]$ and the two time arguments differ by $(1-c_4)(t-t_0)\leq 2S(1-c_4)$, so taking $c_4$ close enough to $1$ makes the bracket at most $c_2/2$ there as well. As in Remark~\ref{rem:average-system-convergence-speed}, this presumes that $\beta$ is the envelope generated from $V$ and $W$ by the Comparison Lemma, which is what licenses the appearance of $c_4$ inside the second argument of $\beta$ in the first place. Now choose
\begin{equation*}
    a^*\left(\alpha_{\beta}^{-1}(c_2/2), \alpha_{\beta}(c_1), c_4\right) = a^*
\end{equation*}
such that, for all $a \in (0,a^*)$,
\begin{equation}
    \label{eq:proof:comparison-functions:vector-field-compact-convergence}
	\left\Vert R_f(\hat{\bar{x}}, a) \right\Vert \left\Vert \frac{\partial V}{\partial \hat{\bar{x}}}(\hat{\bar{x}})\right\Vert \leq (1 - c_4) W(\hat{\bar{x}}) \quad \forall\ \hat{\bar{x}}\in\mathcal{B}_{\alpha_\beta(c_1)}\setminus \mathcal{B}_{\alpha_\beta^{-1}(c_2 /2)}
\end{equation}
This means that 
\begin{equation}
	\frac{\partial V}{\partial \hat{\bar{x}}} \bar{f}(\hat{\bar{x}}, a) \leq -c_4 W(\hat{\bar{x}})\quad \forall\ \hat{\bar{x}}\in\mathcal{B}_{\alpha_\beta(c_1)}\setminus \mathcal{B}_{\alpha_\beta^{-1}(c_2/2)}
\end{equation}
and that the trajectories of the average system satisfy
\begin{align}
	\Vert \hat{\bar{x}}(t_0) \Vert \leq c_1 \implies
		\Vert \hat{\bar{x}}(t) \Vert &\leq
			\max\left\lbrace \beta (\Vert \hat{\bar{x}}(t_0) \Vert, c_4 (t - t_0)),\ \frac{c_2}{2} \right\rbrace \quad \forall\ t\in[t_0,\infty) \\
        \label{eq:average-system-sgpas-conclusion}
		&\leq \beta(\Vert \hat{\bar{x}}(t_0) \Vert, t - t_0) + c_2 \quad \forall\ t\in[t_0,\infty)
\end{align}
where \eqref{eq:average-system-sgpas-conclusion} is a result of $\max\lbrace r_1 , r_2\rbrace \leq r_1 + r_2$ when $r_1,r_2 > 0$ and the inequality in \eqref{eq:beta-closeness}. Thus, we meet the definition of sGPAS w.r.t.~$a$ given in Definition~\ref{def:sgpuas:alternative-definition}.
\end{pf}

We find this proof slightly more opaque than the one based on the Lyapunov function $V$ given in Section~\ref{sec:main-results:average-system} due to the implicit finding of $c_4$ and reliance on $\alpha_\beta$ rather than the bounds $\alpha_1$ and $\alpha_2$ found for $V$. Since the $\mathcal{KL}$ bound for the trajectories is normally estimated from a valid candidate Lyapunov function, one can obtain $\alpha_{\beta} = \alpha_1^{-1}\circ\alpha_2$. However, finding $\beta$ requires solving the differential inequality
\begin{equation}
	\frac{d y}{dt} \leq - \min_{z\in\mathbb{R}^n \text{ s.t. } V(z) = y} W(z)
\end{equation}
or a conservative estimate of it in order to establish an outer bound on all possible $V(t)$ by the Comparison Lemma \cite[Lem 3.4]{ref:khalil-2002}. Solving this differential inequality might be nontrivial, and even if this is accomplished, one must still compute $c_4$ based on the functional constraint given in \eqref{eq:beta-closeness}. Lastly, determining $a^*$ still requires a Lyapunov function for \eqref{eq:proof:comparison-functions:vector-field-compact-convergence} since we cannot use a ``closeness of trajectories'' result due to the model-based system only being continuous and not necessarily locally Lipschitz.

\subsection{The Model-Free System is sGPUAS}

Instead of relying on the $\varepsilon-\delta$ definitions for the proofs in Section~\ref{sec:main-results:model-free}, we provide an alternative proof based on comparison functions, although this proof is still an induction-based one. The outline of the proof is to first establish that there exists $a^*$ and $\omega^*(a)$ such that all trajectories with an initial condition of $\hat{x}(t_0)\in\mathcal{B}_{c_1}$ end up within $\mathcal{B}_{c_5}$ at time intervals $t_1,t_2,\ldots$ where $t_n = t_0 + nT$ for some $T > 0$. The second part is to then establish trajectories with initial conditions $\hat{x}(t_0)\in\mathcal{B}_{c_5}$ are within $\mathcal{B}_{c_2}$ for $t\in[t_1,\infty)$.

\begin{pf}
	Given a $c_1 > 0$ and $c_2>0$, where wlog $c_1 > c_2$, we first choose a $c_5$ such that $\frac{c_5}{3} + \alpha_{\beta}(c_5) < c_2$.	Note that this inequality implies that $c_5 < 3 c_2 / 4$ by $\alpha_{\beta}(r) \coloneqq \beta(r, 0) \geq r$, a condition necessary for $\beta$ to be used in stability definitions. Additionally, since the average system is sGPAS w.r.t.~$a$, there exists an $a^*$ such that, for all $a\in(0,a^*)$,
	\begin{equation}
		\Vert \hat{\bar{x}}(t_0) \Vert \leq c_1 \implies \Vert \hat{\bar{x}}(t) \Vert \leq \beta(\Vert \hat{\bar{x}}(t_0) \Vert, t - t_0) + \frac{c_5}{3}\quad \forall t\in[t_0,\infty)
	\end{equation}
	holds for the average system. Since $\beta\in\mathcal{KL}$, there exists a $T > 0$ such that $\beta(c_1, T) \leq c_5 / 3$ as the second argument of $\beta$ is only a function of the difference $t-t_0$. Given $T$ and
	\begin{equation}
		\Delta = \min\left\lbrace \frac{c_5}{3},\ c_2 - \alpha_{\beta}(c_5) - \frac{c_5}{3}\right\rbrace
	\end{equation} we can establish $\omega^*(T, \Delta, a) = \omega^*$ such that, for all $\omega > \omega^*$, $\Vert \hat{x}(t) - \hat{\bar{x}}(t) \Vert < \Delta$ if $\hat{x}(t_0) = \hat{\bar{x}}(t_0)\in\mathcal{B}_{c_1}$. With this closeness of trajectories, we can establish that, for all $a\in(0,a^*)$ and all $\omega > \omega^*(a)$, 
	\begin{equation}
		\Vert \hat{x}(t_0) \Vert \leq c_1 \implies \Vert \hat{x}(t) \Vert 
			\leq \Vert \hat{\bar{x}}(t) \Vert + \Delta
			\leq \beta(\Vert \hat{x}(t_0) \Vert, t - t_0 ) + \frac{2 c_5}{3}
			\leq \beta(\Vert \hat{x}(t_0) \Vert, t - t_0 ) + c_2 \quad\forall\ t\in[t_0,t_1]
	\end{equation}
	and $\Vert \hat{x}(t_1) \Vert \leq c_5$. Now we consider the time intervals $[t_n,t_{n+1}]$ where $t_n = t_0 + nT$ and $n=1,2,\ldots$. First, we note the trajectory bound
	\begin{equation}
		\Vert \hat{x}(t_n) \Vert \leq c_5 \implies \Vert \hat{x}(t) \Vert 
			\leq \Vert \hat{\bar{x}}(t) \Vert + \Delta
			\leq \beta(c_5, t - t_n ) + \frac{c_5}{3} + \Delta
	\end{equation}
	for $t\in[t_n, t_{n+1}]$. This bound has two implications; the first being $\Vert \hat{x}(t_{n+1}) \Vert \leq c_5$ for $n=1,2,\ldots$ and the second being
	\begin{equation}
		\label{eq:comparison-function-proof:induction-max-bound}
		\Vert \hat{x}(t) \Vert \leq \alpha_{\beta}(c_5) + \frac{c_5}{3} + \Delta \leq c_2
	\end{equation}
	for $t\in[t_n,t_{n+1}]$. However, the combination of the two allows us to extend the time interval for \eqref{eq:comparison-function-proof:induction-max-bound} from a finite time interval to the infinite time interval $t\in[t_1,\infty)$. Combining our analysis on the time intervals $[t_0,t_1]$ and $[t_1,\infty)$, we conclude by induction that, for all $a\in(0,a^*)$ and all $\omega > \omega^*(a)$,
	\begin{equation}
		\Vert \hat{x}(t_0) \Vert \leq c_1 \implies \Vert \hat{x}(t) \Vert \leq \beta(\Vert \hat{x}(t_0) \Vert, t - t_0 ) + c_2\quad \forall\ t\in[t_0,\infty)
	\end{equation}
	which is sGPUAS by Definition~\ref{def:sgpuas:alternative-definition}.
\end{pf}
One key difference between this proof and the one presented in Section~\ref{sec:main-results:model-free} is the conservativeness of $a^*$ and $\omega^*$. Both values must be calculated with respect to some region $\mathcal{D}$ containing the origin in its interior. If we compare two such regions, say $\mathcal{D}_1$ and $\mathcal{D}_2$ with $\mathcal{D}_1 \subset \mathcal{D}_2$, we naturally expect $a^*$ to be smaller and $\omega^*$ to be larger for $\mathcal{D}_2$ than for $\mathcal{D}_1$. Therein lies the difference: proofs with comparison functions are forced to use balls centered at the origin while the proofs in Section~\ref{sec:main-results:average-system} can be altered to utilize the sublevel sets of $V$ instead. For example, the practical boundedness of the average system can be shown with a sublevel of $V$ that is a superset of $\mathcal{B}_{c_1}$. One only needs to find an $a^*$ such that the Lie derivative of $V$ along the average system is negative only at the boundary of the level set. Determining $a^*$ from the sublevel set boundary yields a value that is at least as large as the value obtained using the comparison functions. The balls used in the proof are, at best, the same as the $V$ sublevel sets required for the calculations but often contain regions not within the $V$ sublevel sets. Consequently, the estimates of $a^*$ and $\omega^*$ from the comparison function based approach are generally more conservative than those of the Lyapunov function based approach.

\section*{Declaration of generative AI and AI-assisted technologies in the manuscript preparation process}

During the preparation of this work, the authors used Claude to seek opportunities for expanding upon the authors' exposition of the central ideas. After using this tool/service, the authors reviewed and edited the content as needed and take full responsibility for the content of the published article.

\bibliographystyle{elsarticle-num}
\bibliography{references}


\vskip6pc

\bio{mcnamee-headshot}
	Patrick McNamee received his B.S. ('17) and M.S. ('19) in Aerospace Engineering and M.S. ('21) in Computer Science from the University of Kansas. He is currently a graduate research assistant with the Dynamical Systems and Intelligent Laboratory at San Diego State University where he is pursuing a Ph.D. in Mechanical Engineering. His research interests include autonomous vehicles and extremum seeking control.
\endbio

\vskip3pc

\bio{Zahra-Nili-Ahmadabadi}
	Zahra Nili Ahmadabadi is Assistant Professor with the Mechanical Engineering Department at San Diego State University (SDSU). She is a recipient of the ARO Early career award. Nili is currently Associate Editor of Sage Journal of Vibration and Control. Her current research interests include robot perception and control, and nonlinear systems.
\endbio

\vskip6pc

\bio{Krstic-red-tie}
	Miroslav Krsti{\'c} is Distinguished Professor and Senior Associate Vice Chancellor for Research at University of California, San Diego. He is Fellow of IEEE, IFAC, ASME, SIAM, AAAS, AIAA (Assoc. Fellow), IET (UK), and member of the Serbian Academy of Sciences and Arts. Krsti{\'c} is a recipient of the Bellman Control Heritage Award, Bode Lecture Prize, SIAM Reid Prize, ASME Oldenburger Medal, Nyquist Lecture Prize, Paynter Award, Ragazzini Education Award, IFAC Nonlinear Control Systems Award, IFAC Ruth Curtain Distributed Parameter Systems Award, IFAC Adaptive and Learning Systems Award, Chestnut textbook prize, CSS Distinguished Member Award, the PECASE, NSF Career, and ONR Young Investigator awards, and the Schuck (’96 and ’19) and Axelby paper prizes. Krsti{\'c} is also the inaugural recipient of the A.V. Balakrishnan Award for the Mathematics of Systems. He is Editor-in-Chief of IEEE Transactions on Automatic Control and has coauthored 19 books and about 500 journal papers.
\endbio

\end{document}